# CONDITIONAL PREDICTIVE INFERENCE POST MODEL SELECTION

By Hannes Leeb

*Yale University*

We give a finite-sample analysis of predictive inference procedures after model selection in regression with random design. The analysis is focused on a statistically challenging scenario where the number of potentially important explanatory variables can be infinite, where no regularity conditions are imposed on unknown parameters, where the number of explanatory variables in a "good" model can be of the same order as sample size and where the number of candidate models can be of larger order than sample size. The performance of inference procedures is evaluated conditional on the training sample. Under weak conditions on only the number of candidate models and on their complexity, and uniformly over all data-generating processes under consideration, we show that a certain prediction interval is approximately valid and short with high probability in finite samples, in the sense that its actual coverage probability is close to the nominal one and in the sense that its length is close to the length of an infeasible interval that is constructed by actually knowing the "best" candidate model. Similar results are shown to hold for predictive inference procedures other than prediction intervals like, for example, tests of whether a future response will lie above or below a given threshold.

## 1. Introduction.

1.1. *Motivation and summary.* This paper is about inference on future observations based on a model that has been selected on the basis of the data and then fitted to the same data. We focus, in particular, on situations where the number of candidate models is large and where the number of explanatory variables in a "good" model can be large as well, in relation to sample size. Such a situation is faced, for example, by Stenbakken and Souders [31] who predict the performance of analog/digital converters from partial measurements by selecting 64 explanatory variables (measurements) from a total of 8192 based on a sample of size 88; further examples include [1, 8, 12, 30, 33, 34, 35] and [37]. Note that, in these studies, the model that is selected, on the basis of the data, is often quite complex in relation to sample size, in the sense that the number of explanatory variables in the selected model and the sample size are of the same order of magnitude. Also note that the total number of candidate models in these studies exceeds sample size by several orders of magnitude. In such situations, inferential tools that assess the predictors' accuracy like, for example, the mean-squared error of the predictor, or prediction intervals, are needed.

We consider a Gaussian regression model with random design, where the number of explanatory variables can be infinite, and where no regularity conditions are imposed on the unknown parameters. We use a variant of generalized cross-validation to evaluate the performance of candidate models for prediction out-of-sample,[1] to select a "good" model and to conduct predictive inference based on the selected model. The performance of the resulting model selector and the quality of predictive inference procedures are evaluated conditional on the training sample. We describe the performance of these methods by explicit finite-sample performance bounds. For example, we show that the proposed prediction interval is approximately valid and short, with high probability, even in statistically challenging situations where the number of explanatory variables in a "good" model is of the same order as sample size, and where the total number of candidate models is of a larger order than sample size. Here, approximately valid means that the prediction interval's actual coverage probability is close to the nominal one, and approximately short means that its length is close to the length of a certain infeasible "prediction interval" that is based on actually knowing the "best" candidate model. Our results hold uniformly over all data-generating processes under consideration.

1.2. *Our results in broader context.* In the literature, results on predictive inference after model selection are scarce.[2] The finite-sample distribution of a linear predictor based on the selected model can be computed explicitly in sufficiently simple settings (see Leeb [18] and [19]). However, these results only allow for rather restricted collections of candidate models; moreover, as the number of candidate models increases, the resulting formulae get increasingly complicated and computationally infeasible. From the perspective of traditional large-sample analyses, on the other hand, predictive inference after model selection is typically rather trivial. Consider, for example, a parametric linear model where the response is a linear function of a finite number of explanatory variables and a random disturbance. Under standard assumptions, every sensible model selection procedure typically leads to a post model selection estimator that is consistent, even uniformly consistent (see, e.g., Propositions A.9 and B.1 in [22]). In large samples, the random disturbance is therefore the dominant source of error when predicting a new response. Thus, as far as prediction is concerned, all sensible model selection procedures perform alike in parametric settings in the large-sample limit. The same is true in nonparametric settings for appropriately chosen estimators of the true regression function, provided that the regression function is a priori restricted to a sufficiently regular family like, say, a Besov body or a collection thereof, as it is often considered in nonparametric function estimation. In situations where the true regression parameter or function can be estimated consistently or uniformly consistently, research is typically focused on finding estimators with good convergence rates, or on finding confidence sets for the true regression parameter or function that are valid and small.

---

[1]Here, prediction "out-of-sample" means prediction of new responses given hitherto unobserved explanatory variables, whereas "in-sample" prediction means prediction of new responses for the same explanatory variables as observed in the training data.

[2]This is in spite of the fact that predictive inference by itself is a rather well researched field (see, e.g., [3] and [9] for a frequentist and a Bayesian approach, respectively, as well as the references given there).

In this paper, we consider predictive inference after model selection in a situation that is difficult to analyze by exact finite-sample results or by large-sample limit theory. In particular, we focus on the statistically challenging scenario where the number of explanatory variables in a "good" model can be of the same order as sample size, and where the number of candidate models can be of larger order than sample size. This situation is typically too complex for an exact finite-sample analysis. Also, this situation is such that large-sample limit approximations cannot be guaranteed to be accurate in most cases. We do not rule out the case where: (i) a very simple model fits well and (ii) the number of candidate models is small, in relation to sample size. However, our results are most interesting in the case where one of these two conditions is not met.

There are, however, a couple of results, in both parametric and nonparametric settings, that indicate that inference after model selection is a hard problem that is subject to certain insurmountable obstructions. Most of these results consider inference on the regression parameter itself, on components thereof or on the mean of a future response.

Consider, first, a parametric linear model, with Gaussian errors and fixed design (under standard assumptions), and a linear predictor that is constructed based on the outcome of a data-driven model selection step. It is well known that the distribution of such a linear predictor, properly scaled and centered, typically depends on unknown parameters in a nontrivial way and can be highly nonnormal, regardless of sample size (see [21, 27]).[3] Moreover, Leeb and Pötscher [23, 24] showed that the distribution of such a linear predictor cannot be estimated in a uniformly consistent fashion, except in degenerate and trivial cases. Concerning confidence intervals for the mean of a future response, the results of Joshi [15] entail in the known-variance case that the standard interval based on fitting the overall model is admissible and uniquely minimax with respect to a loss function that measures both coverage probability and interval size (in the class of all randomized Lebesgue measurable confidence sets and up to trivial equivalences). (See also [17] for further references and results on confidence intervals post model selection.)

In nonparametric function estimation, there are well-known limits to the adaptivity of honest $L_2$ confidence balls for the true regression function. Here, "honest" means that the confidence ball guarantees coverage probability over the whole function space under consideration, and "adaptivity" means that smoother regression functions (i.e., functions that belong to restricted submodels) are covered by smaller balls. In essence, larger function spaces limit the amount of adaptivity possible. This was first discovered by Li [25] and was further analyzed in [2, 4, 6, 7, 10, 14, 16] and [28]. Moreover, Baraud [2] also shows that honest and short confidence balls are feasible only if the error variance is assumed to lie in some bounded subset of $(0, \infty)$, and that loose variance bounds close to zero or infinity lead to large confidence balls. If an honest confidence band (i.e., an $L_\infty$ confidence ball) is desired, then the limits to adaptivity are even more pronounced (see [11]).

In the setting of this paper, where the goal is prediction out-of-sample, we demonstrate that prediction intervals post model selection can be simultaneously valid and short in an approximate sense and with high probability, irrespective of unknown parameters. The proposed prediction interval has the following two properties, except on an event whose probability is bounded by the expression in (1), which follows: (i) Its actual coverage probability is close to the nominal coverage probability. (ii) Its length is close to the length of a certain infeasible shortest possible interval that is constructed from actually knowing the "best" candidate model. These statements hold uniformly over all data-generating processes under consideration.

On a technical level, this paper is related to Breiman and Freedman [5] in two regards: First, the model considered in this paper contains the model considered in [5] as a special case [see the discussion following (2)], and, second, our results rely on a corresponding extension of Theorem 1.3 of [5] (see Proposition 2.1 and the attending discussion). The results derived here, however, differ considerably from those of [5] in terms of scope and content. We allow for families of candidate models of essentially arbitrary size and structure, while [5] is focused on up to $n/2$ models that are nested (where $n$ denotes sample size). Moreover, we give finite-sample results that hold uniformly over all data-generating processes under consideration, while the main result in [5] is a pointwise large-sample limit result that requires that the true regression parameter has infinitely many nonzero components.

---

[3]This fact is at odds with a result of Shen, Huang and Ye [29], which claims that the limit distribution of a post model selection estimator in a parametric setting is normal with mean zero and estimable variance/covariance matrix (see Theorems 3 and 4 in that paper). Inspection of that paper reveals that the proof of Theorem 3 is in error as it stands, and that said theorem does not hold as claimed. Private communication with the authors has confirmed this.

1.3. *Outline of the paper.* As the data-generating process, we consider a Gaussian linear model with random design that is described in Section 2, where the number of potentially important explanatory variables can be infinite. We assume that the error and also the explanatory variables are jointly Gaussian, like Breiman and Freedman [5]. Assuming the data to be Gaussian allows us to derive explicit finite-sample performance bounds by relatively elementary means and to clearly showcase the mechanisms underlying our results. Simulation results in [20] strongly suggest that the assumption of Gaussianity is not essential, and unpublished preliminary results, which rely on random matrix theory, point in the same direction. The unknown parameters in this setting are the sequence of regression coefficients as well as the means and the variance/covariance structure of the explanatory variables and of the error term. No additional regularity conditions are imposed on the unknown parameters.

We consider a scenario where the model is selected and fitted to the data once and is then used repeatedly for prediction and for predictive inference. For performance measures, like the mean-squared error of a predictor or the coverage probability or length of a prediction interval, we therefore adopt a conditional perspective and treat the training sample as fixed and the future response and its corresponding explanatory variables as random.[4]

Given a sample of size $n$ and a collection $\mathcal{M}$ of candidate models, a preliminary first goal is to evaluate models $m \in \mathcal{M}$ based on their performance for prediction out-of-sample, and to select a model that performs well for this purpose; this is the subject of Section 3. Our second and main goal is to conduct inference on future observations based on the selected model like, for example, prediction intervals; this goal is studied in Section 4.

---

[4]This deviates from conventional linear model theory, where, usually, the training sample is considered random and where, often, the explanatory variables that are used for prediction are considered as fixed. Regarding prediction intervals, our approach may be compared to the average coverage probability introduced by Wahba [36] and further analyzed in [26].

To achieve both goals outlined in the preceding paragraph, we consider a model selector and predictive inference procedures post model selection that are based on a variant of generalized cross-validation (and that are described in detail later). We show that the proposed prediction interval is approximately valid and short, except on an event whose probability is bounded by

$$(1) \qquad C_1 \exp[\log \#\mathcal{M} - C_2(n - |\mathcal{M}|)],$$

uniformly over all data-generating processes under consideration. Here, $\#\mathcal{M}$ is the number of candidate models, $|\mathcal{M}|$ is the number of explanatory variables in the most complex candidate model and $C_1$ and $C_2$ are explicit positive constants. The bound in (1) decreases exponentially fast in $n - |\mathcal{M}|$ and increases only linearly in $\#\mathcal{M}$. This allows for very large classes of potentially complex candidate models. If the upper bound in (1) is small, the proposed prediction interval is approximately valid and its length is close to that of a certain infeasible "prediction interval" that is based on actually knowing the "best" candidate model, with high probability (see Propositions 4.3 and 4.4 on page 15 for details). Furthermore, we show that the following statements hold, except on an event whose probability is bounded by (1) (with different values of the constants $C_1$ and $C_2$): (i) The performance of the selected model is close to the performance of the "best" candidate model; (ii) the estimated performance of the selected model is close to its actual performance; and (iii) in general, the proposed procedures for predictive inference post model selection are approximately valid. In a simulation example, we use a training sample of 2000 observations to perform a greedy search through a pool of over $10^{15}$ candidate models (see Section 5).

**2. Basic assumptions and quantities of interest.** For the data-generating process, consider a response $y$ that is related to a collection of explanatory variables $(x_j)_{j \geq 1}$ by

$$(2) \qquad y = \sum_{j=1}^{\infty} x_j \beta_j + u.$$

Assume that the model includes an intercept (i.e., $x_1 = 1$) and that the $x_j$'s for $j > 1$ and $u$ are jointly nondegenerate Gaussian with unknown means and variance/covariance structure, such that the sum converges in $L_2$.[5] No additional regularity conditions will be imposed on the data-generating process throughout the paper. Breiman and Freedman [5] consider a special case of (2), where the mean of the explanatory variables is known (and equal to zero), and where no intercept is included (i.e., $\beta_1 = 0$).

---

[5]Hence, the distribution of any finite subset of $\{x_j : j > 1\} \cup \{u\}$ is a nondegenerate Gaussian with unknown mean-vector and variance/covariance matrix. It is often also assumed that the $x_j$'s are uncorrelated with $u$ and that $u$ has mean zero. This assumption is not needed here. In essence, $u$ plays the role of as an unobserved explanatory variable.

The minimal requirement, that the right-hand side of (2) converges in squared mean, restricts the possible values of $\beta = (\beta_j)_{j \geq 1}$ in a way that depends on the moments of the explanatory variables. For example, if the $x_j$'s, $j > 1$, are independent and identically distributed with mean zero and variance, say, one, then $\beta$ can be any sequence of coefficients in $l_2$. This shows that (2) covers a large class of data-generating processes; further examples are outlined in Remark 6.1. Of course, (2) also covers parametric models with only finitely many explanatory variables (i.e., the case where $\beta_j = 0$ from some index onward). Moreover, the requirement of nondegeneracy can be relaxed as outlined in Remark 6.2.

Consider a sample of size $n$ from (2). The sample will be denoted by $(Y, X)$ with $Y$ denoting the $n$-vector $Y = (y^{(1)}, \ldots, y^{(n)})'$ and $X$ denoting the $n \times \infty$ "matrix" or net $X = (x^{(1)}, \ldots, x^{(n)})'$, where $(y^{(i)}, x^{(i)})$ are independent and identically distributed (i.i.d.) copies of $(y, x)$ as in (2).

The training sample will be used to fit finite-dimensional submodels of (2) that restrict some coefficients of $\beta$ to zero, where the intercept $\beta_1$ is always left unrestricted. Each such submodel is described by a 0–1 sequence $m = (m_j)_{j \geq 1}$, where $m_j = 0$ if the $j$th coefficient of $\beta$ is restricted to zero and $m_j = 1$, otherwise. The number of unrestricted regression coefficients (i.e., $\sum_{j \geq 1} m_j$) is denoted by $|m|$. We assume that $|m| < n - 1$ throughout the paper.

Consider a finite collection of candidate models that will be denoted by $\mathcal{M}$ (which, of course, may depend on sample size $n$). Assume that each model $m \in \mathcal{M}$ satisfies

(3) $$m_1 = 1 \quad \text{and} \quad |m| < n - 1$$

as before.[6] We write $|\mathcal{M}|$ for the number of parameters in the most complex model in $\mathcal{M}$; that is,

$$|\mathcal{M}| = \max_{m \in \mathcal{M}} |m|,$$

and we write $\#\mathcal{M}$ for the number of candidate models in $\mathcal{M}$.

For later use, let $\sigma^2(m)$ denote the variance of $y$ conditional on those explanatory variables that are included in the model $m$; that is,

$$\sigma^2(m) = \text{Var}[y|x_j : m_j = 1, j \geq 1]$$

---

[6]In practice, the choice of candidate models $\mathcal{M}$ to consider at sample size $n$ is often guided by prior knowledge or suspicions about the structure of the underlying parameters. For example, if it is assumed or suspected that the coefficients of $\beta$ are sparse in an appropriate sense, one might consider appropriately sparse candidate models a well; such a case is discussed in the simulation example in Section 5. Another example is the case where the coefficients of $\beta$ are assumed or suspected to taper off at a certain rate.

for each $m \in \mathcal{M}$. Note that this conditional variance does not depend on the $x_j$'s because $y$ and $(x_j)_{j\geq 1}$ are jointly Gaussian. Also note that $0 < \operatorname{Var}[u] \leq \sigma^2(m) \leq \operatorname{Var}[y]$ for each data-generating process as in (2); in particular, $\sigma^2(m)$ is always positive.

The least-squares method will be used to fit models to the training sample.[7] The restricted least-squares estimator corresponding to a model $m \in \mathcal{M}$ is denoted by $\hat{\beta}(m) = (\hat{\beta}_j(m))_{j \geq 1}$ and is defined as follows: For $j$ satisfying $m_j = 0$, $\hat{\beta}_j(m)$ equals zero; the $|m|$ remaining components of $\hat{\beta}(m)$ are obtained by regressing $Y$ on the observed values of those regressors that are included in the model $m$ (on the probability zero event where the resulting $n \times |m|$ regressor matrix is rank deficient, we use the Moore–Penrose inverse, say, in the least-squares formula). The usual variance estimator based on model $m$ will be denoted by $\hat{\sigma}^2(m)$ and is given by $\hat{\sigma}^2(m) = (n - |m|)^{-1} \operatorname{RSS}(m)$ with $\operatorname{RSS}(m)$ denoting the residual sum of squares obtained by fitting model $m$ to the training sample [note that $\hat{\sigma}^2(m) > 0$, almost surely].

The performance of a model will be evaluated in terms of the conditional mean-squared error of the linear predictor obtained from fitting the model to the training sample. Let $(y^{(f)}, x^{(f)})$ be a new copy of $(y, x)$ as in (2), independent of the sample $(Y, X)$. Based on a model $m \in \mathcal{M}$ and the sample $(Y, X)$, the usual least-squares predictor of $y^{(f)}$ will be denoted by $\hat{y}^{(f)}(m)$ and is given by

$$\hat{y}^{(f)}(m) = \sum_{j=1}^{\infty} x_j^{(f)} \hat{\beta}_j(m).$$

Note that all but $|m|$ coefficients of the restricted least-squares estimator $\hat{\beta}(m)$ are zero. The conditional mean-squared error of the predictor is now defined as

$$\rho^2(m) = E[(\hat{y}^{(f)}(m) - y^{(f)})^2 | Y, X].$$

Note that $\rho^2(m)$ depends on the training sample and, hence, is a random variable, and that $\rho^2(m)$ also depends on the unknown parameters in (2). In particular, $\rho^2(m)$ is unknown. We will also consider the corresponding unconditional mean-squared error of the predictor (i.e., $E[\rho^2(m)]$) and the positive square root $\rho(m)$ of $\rho^2(m)$.

If the predictor $\hat{y}^{(f)}(m)$ is to be used for inferences about a new response $y^{(f)}$, the distribution of the prediction error $\hat{y}^{(f)}(m) - y^{(f)}$ is of particular interest. Conditional on the training sample, the distribution of the prediction error $\hat{y}^{(f)}(m) - y^{(f)}$ will be denoted by $\mathbb{L}(m)$. Clearly, $\mathbb{L}(m)$ is a Gaussian, and we write $\nu(m)$ for the mean of that distribution and $\delta(m)$ for its standard deviation. In other words,

$$\hat{y}^{(f)}(m) - y^{(f)}|Y, X \sim N(\nu(m), \delta^2(m)) \equiv \mathbb{L}(m).$$

---

[7]While it is tempting to also consider penalized least-squares or more general shrinkage estimators, particularly for complex candidate models, our current methods cannot handle these estimators.

Note that $\nu(m)$ is also the conditional bias of the predictor $\hat{y}^{(f)}(m)$, conditional on the training sample. As before, also note that the distribution of these quantities depends on the unknown parameters in (2), so that $\nu(m)$ and $\delta^2(m)$ are unknown. Of course, we have $\rho^2(m) = \nu^2(m) + \delta^2(m)$.

In terms of the conditional mean-squared error of prediction, the best candidate model is a minimizer of $\rho^2(m)$ over $m \in \mathcal{M}$. We write $m_\rho$ for such a minimizer; that is,

$$m_\rho = \underset{m \in \mathcal{M}}{\arg\min}\, \rho^2(m)$$

(on the event of multiple minimizers, $m_\rho$ is taken as a measurable selection from the set of minimizers). In Section 4, we will also consider the candidate model for which the conditional distribution of the prediction error [i.e., $\mathbb{L}(m)$] is most concentrated. That model [i.e., a measurable minimizer of $\delta^2(m)$ over $m \in \mathcal{M}$] is denoted by $m_\delta$.

For deriving and analyzing estimators for the quantities of interest $\delta^2(m)$, $\nu^2(m)$ and $\rho^2(m)$ and for understanding the mechanisms underlying our main findings, the following result will be instrumental.

PROPOSITION 2.1. *For each fixed model $m \in \mathcal{M}$, the conditional variance of the prediction error $\hat{y}^{(f)}(m) - y^{(f)}$ given $(Y, X)$ [i.e., $\delta^2(m)$] has the same distribution as $\sigma^2(m)$ multiplied by the sum of one and the ratio of two independent chi-square random variables with $|m| - 1$ and $n - |m| + 1$ degrees of freedom, respectively:*

$$\delta^2(m) \sim \sigma^2(m)\left(1 + \frac{\chi^2_{|m|-1}}{\chi^2_{n-|m|+1}}\right).$$

*The conditional bias of $\hat{y}^{(f)}(m)$ given $(Y, X)$ has mean zero (i.e., $E[\nu(m)] = 0$). Moreover, the squared conditional bias $\nu^2(m)$ has the same distribution as $\delta^2(m)/n$ multiplied by an independent chi-square random variable with one degree of freedom:*

$$\nu^2(m) \sim \frac{\chi^2_1}{n}\sigma^2(m)\left(1 + \frac{\chi^2_{|m|-1}}{\chi^2_{n-|m|+1}}\right).$$

*Finally, the usual variance estimator in model $m$ is distributed as $\hat{\sigma}^2(m) \sim \sigma^2(m)\chi^2_{n-|m|}/(n - |m|)$ [in case $|m| = 1$, the expression $\chi^2_{|m|-1}$ in the preceding two displays is to be interpreted as constant equal to zero, so that $\delta^2(m) = \sigma^2(m)$ and $\nu^2(m) \sim (\chi^2_1/n)\sigma^2(m)$ in this case].*

Proposition 2.1 extends Theorem 1.3 of Breiman and Freedman [5], which describes the distribution of $\delta^2(m)$ in the case where the regressors in (2) all have mean zero and where models do not contain an intercept. For a slightly different conditioning sigma-field, the distribution of the corresponding conditional mean-squared error of the predictor is also derived by Thompson [32].

Proposition 2.1 shows that the squared conditional bias $\nu^2(m)$ is of smaller order, by a factor of $1/n$ and in probability, than the conditional variance $\delta^2(m)$. A little reflection shows that this is no surprise, for example, in the case where the fitted model is correct (i.e., in the case where $m$ contains all nonzero coefficients of $\beta$). By Proposition 2.1, the same is true regardless of how well the fitted model describes the true one. Of course, $\nu^2(m)$ can be substantial because of either overfit or underfit, say. But, irrespective of that, the conditional variance $\delta^2(m)$ is the dominating factor in $\rho^2(m) = \nu^2(m) + \delta^2(m)$, in probability. Another feature revealed by Proposition 2.1 is that the distributions of $\nu^2(m)$ and $\delta^2(m)$ depend on the unknown parameters in (2) only through $\sigma^2(m)$, and that $\sigma^2(m)$ can be estimated from the training sample with good accuracy, provided only that $n - |m|$ is large. For later use, we can also read-off the expected values of $\delta^2(m)$, $\nu^2(m)$ and $\rho^2(m)$ from Proposition 2.1. Because the mean of $1/\chi^2_{n-|m|+1}$ equals $1/(n-|m|-1)$ for $n - |m| - 1 > 0$, the mean of $\delta^2(m)$ equals $\sigma^2(m)(n-2)/(n-|m|-1)$ and the mean of $\nu^2(m)$ is $n^{-1}\sigma^2(m)(n-2)/(n-|m|-1)$. From this, we also see that the mean of $\rho^2(m)$, that is, the (unconditional) mean-squared error of the predictor $\hat{y}^{(f)}(m)$, is given by

$$E[\rho^2(m)] = \sigma^2(m)\frac{n-2}{n-1-|m|}\left(1 + \frac{1}{n}\right).$$

This formula for $E[\rho^2(m)]$ is also derived in [13] and [32] by different means. Finally, Proposition 2.1 suggests that $\delta^2(m)/E[\delta^2(m)]$ is close to one provided only that $n - |m|$ is sufficiently large, and that the same is true for $\rho^2(m)/E[\rho^2(m)]$. Formalizing this idea and using variations of Chernoff's method will lead to the main results of this paper, Theorems 3.1 and 4.1, which follow.

**3. Evaluating and selecting models.** The performance of model $m$, as measured by the conditional mean-squared error of the predictor $\hat{y}^{(f)}(m)$ [i.e., as measured by $\rho^2(m)$] depends on unknown parameters and, hence, cannot be used directly for model selection. We now consider several estimators for $\rho^2(m)$. In view of Proposition 2.1 and the ensuing formula for $E[\rho^2(m)]$, we see that an unbiased estimator for $E[\rho^2(m)]$ is given by

$$\check{\rho}^2(m) = \hat{\sigma}^2(m)\frac{n-2}{n-1-|m|}\left(1 + \frac{1}{n}\right)$$

(see also [13, 32]). Of course, this estimator is also unbiased for $\rho^2(m)$. The estimator $\breve{\rho}^2(m)$ is closely related to two well-known model selectors, namely generalized cross-validation and the $S_p$ criterion, whose objective functions are defined by

$$\mathrm{GCV}(m) = \hat{\sigma}^2(m) \frac{n}{n-|m|} \quad \text{and} \quad S_p(m) = \hat{\sigma}^2(m) \frac{n-2}{n-1-|m|},$$

respectively. For fixed sample size $n$, choosing a model $m$ that minimizes $\breve{\rho}^2(m)$ is equivalent to choosing a model that minimizes $S_p(m)$. Moreover, for most practical purposes, the difference between $\breve{\rho}^2(m)$, $S_p(m)$ and $\mathrm{GCV}(m)$ will be negligible. Because of technical reasons, we consider another estimator that is closely related to the three discussed so far. That estimator will be denoted by $\hat{\rho}^2(m)$ and is given by

$$\hat{\rho}^2(m) = \hat{\sigma}^2(m) \frac{n}{n+1-|m|}.$$

Again, note that the difference between $\hat{\rho}^2(m)$, $\breve{\rho}^2(m)$, $\mathrm{GCV}(m)$ and $S_p(m)$ will be negligible for most practical purposes. The relation between $\mathrm{GCV}(m)$ or $S_p(m)$ and other well-known model selection criteria in our setting is discussed in detail in Section 3.3 of [20]. The next result describes the performance of $\hat{\rho}^2(m)$ as an estimator for the conditional mean-squared error of the predictor [i.e., as an estimator for $\rho^2(m)$] in finite samples.[8]

THEOREM 3.1.  *Fix a candidate model $m \in \mathcal{M}$. For each $\varepsilon > 0$, we have*

(4)  $$P\left(\left|\log \frac{\hat{\rho}^2(m)}{\rho^2(m)}\right| > \varepsilon\right) \leq 6 \exp\left[-\frac{n-|m|}{8} \frac{\varepsilon^2}{\varepsilon + 8}\right].$$

*The relation in the preceding display holds uniformly over the set of all data-generating processes as in (2).*

Theorem 3.1 shows that the estimated performance of model $m$ is close to its true performance, in the sense that the ratio $\hat{\rho}^2(m)/\rho^2(m)$ is close to one with high probability, provided only that $n - |m|$ is large enough, independently of the unknown parameters. The theorem places no restriction on sample size $n$ and on the candidate model $m$ [except for (3) that is maintained throughout the paper]. However, the result is most interesting in the case where the sample size is relatively small compared to the number of parameters in the model, in the sense that $|m|/n$ is not close to zero.

---

[8]In Theorem 3.1, the expression $|\log \hat{\rho}^2(m)/\rho^2(m)|$ is, of course, well defined in case $\hat{\rho}^2(m) > 0$ or, equivalently, in case $\hat{\sigma}^2(m) > 0$, which is an almost sure event. In case $\hat{\rho}^2(m) = 0$, $|\log \hat{\rho}^2(m)/\rho^2(m)|$ is to be interpreted as $\infty$. The same convention is also used, mutatis mutandis, in the results that follow.

In that case, other model selectors like, say, AIC, AICc, FPE or BIC, give a distorted picture of the model's performance, and the model selected by one of these model selection criteria can be anything from mildly suboptimal to completely unreasonable, depending on unknown parameters. These phenomena are discussed at length in Section 3.3 of [20] for the special case where the regressors in (2) are centered to have mean zero and where candidate models do not include an intercept. That discussion also applies to the setting that is considered here, mutatis mutandis.

Because the upper bound in (4) decreases exponentially fast in $n - |m|$, Theorem 3.1 can be used together with Bonferroni's inequality to describe the performance of $\hat{\rho}^2(m)$ when this estimator is used to evaluate the performance of several candidate models. For the collection $\mathcal{M}$ of candidate models introduced at the end of Section 2, recall that model $m_\rho$ minimizes $\rho^2(m)$ over $m \in \mathcal{M}$. The truly best model $m_\rho$ is of course infeasible, but Theorem 3.1 suggests that $\hat{\rho}^2(m)$ can be taken as a proxy for $\rho^2(m)$. Define the empirically best model $\hat{m}$ as a (measurable) minimizer of $\hat{\rho}^2(m)$ over $\mathcal{M}$; that is,

$$\hat{m} = \arg\min_{m \in \mathcal{M}} \hat{\rho}^2(m).$$

For the next result, recall that $|\mathcal{M}|$ denotes the number of parameters in the most complex candidate model and that $\#\mathcal{M}$ denotes the total number of candidate models.

COROLLARY 3.2. *For each $\varepsilon > 0$ and uniformly over all data-generating processes as in (2), we have*

$$(5) \qquad P\left(\log \frac{\rho^2(\hat{m})}{\rho^2(m_\rho)} > \varepsilon\right) \leq 6 \exp\left[\log \#\mathcal{M} - \frac{n - |\mathcal{M}|}{16} \frac{\varepsilon^2}{\varepsilon + 16}\right]$$

*and*

$$(6) \qquad P\left(\left|\log \frac{\hat{\rho}^2(\hat{m})}{\rho^2(\hat{m})}\right| > \varepsilon\right) \leq 6 \exp\left[\log \#\mathcal{M} - \frac{n - |\mathcal{M}|}{8} \frac{\varepsilon^2}{\varepsilon + 8}\right].$$

The first inequality of Corollary 3.2 relates the performance of the empirically best model (i.e., $\hat{m}$) to that of the actually best candidate model (i.e., $m_\rho$) in terms of the relative performance $\rho^2(\hat{m})/\rho^2(m_\rho)$; if the upper bound in (5) is small, the performance of $\hat{m}$ is close to that of $m_\rho$ with high probability. In that case, one can select a "good" model on the basis of the data with high probability. Moreover, the second inequality shows that the performance of the selected model can be estimated accurately, in terms of the relative error $|\log \hat{\rho}^2(\hat{m})/\rho^2(\hat{m})|$, with high probability, provided that the upper bound in (6) is small. It should be noted that the upper bounds in

(5) and (6) do not depend on unknown parameters but only on sample size, on the number of candidate models, and on the number of parameters in the most complex candidate model (i.e., on $n$, $\#\mathcal{M}$ and $|\mathcal{M}|$). In particular, these upper bounds are small if the degrees of freedom in the most complex candidate model (i.e., $n - |\mathcal{M}|$) is sufficiently large compared to $\log \#\mathcal{M}$. This allows for very large classes of potentially very complex candidate models [see also Remark 6.3 for a discussion of the role of the constants $\#\mathcal{M}$ and $|\mathcal{M}|$ in the upper bounds (5), (6), and in the results that follow]. Finally, we note that we actually establish a slightly stronger result during the proof of Corollary 3.2, namely that the result continues to hold with the left-hand side of (6) replaced by $P(\max_{m \in \mathcal{M}} |\log \hat{\rho}^2(m)/\rho^2(m)| > \varepsilon)$. In other words, if the upper bound in (6) is small, then $\hat{\rho}^2(m)/\rho^2(m)$ is close to one for each $m \in \mathcal{M}$ with high probability. In that case, $\hat{\rho}^2(\cdot)$ can be used to approximate the predictive performance not only of $\hat{m}$ but also of other model selection procedures that differ from $\hat{m}$ (see [20] for some examples and further discussion).

The results presented so far are concerned with relative errors like, for example, $\log \hat{\rho}^2(m)/\rho^2(m)$. Theorem 3.1 also entails similar results for absolute errors like, for example, $\hat{\rho}^2(m) - \rho^2(m)$, that parallel results in [20] and are omitted here for the sake of brevity.

**4. Predictive inference based on the selected model.** To use the predictor $\hat{y}^{(f)}(m)$ for inferences about the unseen future response $y^{(f)}$, like prediction intervals for example, the distribution of the prediction error [i.e., of $\hat{y}^{(f)}(m) - y^{(f)}$] is an object of particular interest. Recall that we write $\mathbb{L}(m)$ for the conditional distribution of this prediction error given the training sample. For a fixed candidate model $m$ and fixed training sample, $\mathbb{L}(m)$, of course, depends on unknown parameters and, hence, needs to be estimated; in particular, we need to estimate the conditional bias and the conditional variance of the predictor. Proposition 2.1 shows that (unconditionally) unbiased estimators of $\nu(m)$ and of $\delta^2(m)$ are given by zero and by

$$\check{\delta}^2(m) = \hat{\sigma}^2(m) \frac{n-2}{n-1-|m|},$$

respectively (i.e., $E[\nu(m)] = E[\check{\delta}^2(m) - \delta^2(m)] = 0$). This suggests that the distribution in question [i.e., $\mathbb{L}(m) \equiv N(\nu(m), \delta^2(m))$] might be estimated by $N(0, \check{\delta}^2(m))$. For technical reasons, we consider a slightly different estimator. In particular, we estimate $\delta^2(m)$ by

$$\hat{\delta}^2(m) = \hat{\sigma}^2(m) \frac{n}{n+1-|m|}$$

[which coincides with the estimator $\hat{\rho}^2(m)$ discussed in Section 3], and we estimate the conditional distribution of the prediction error [i.e., $\mathbb{L}(m)$] by

$$\hat{\mathbb{L}}(m) \equiv N(0, \hat{\delta}^2(m)).$$

The next result describes the finite-sample performance of $\hat{\mathbb{L}}(m)$ as an estimator for $\mathbb{L}(m)$ in terms of the total variation distance.

THEOREM 4.1. *Fix a candidate model $m \in \mathcal{M}$. For the conditional distribution of the prediction error of the predictor $\hat{y}^{(f)}(m)$, conditional on the training sample, and for its estimated version [i.e., for $\mathbb{L}(m)$ and for $\hat{\mathbb{L}}(m)$] we have*

$$(7) \qquad P\bigg(\|\hat{\mathbb{L}}(m) - \mathbb{L}(m)\|_{\text{TV}} > \frac{1}{\sqrt{n}} + \varepsilon\bigg) \le 7\exp\bigg[-\frac{n-|m|}{2}\frac{\varepsilon^2}{\varepsilon+2}\bigg]$$

*for each $\varepsilon$ with $0 < \varepsilon \le \log(2)$. The upper bound in the preceding display holds uniformly over the set of all data-generating processes as in (2).*

REMARK 4.1. Because the total variation distance of two probability measures is at most 1, the condition that $\varepsilon$ is at most $\log(2) \approx 0.69$ maintained by Theorem 4.1 is rather innocuous. Inspection of the proof of Theorem 4.1 shows that one can obtain a slightly improved upper bound that also holds for all $\varepsilon > 0$. The downside of this is that the improved upper bound is much more complicated and less revealing.

By Theorem 4.1, the estimated distribution $\hat{\mathbb{L}}(m)$ is close to the true distribution $\mathbb{L}(m)$ in total variation with high probability, provided only that $n - |m|$ is large enough, independently of the unknown parameters. While the theorem places no restrictions on sample size and on the candidate model $m$ [except for (3)], the result is most interesting in the case where the candidate model is relatively complex in the sense that $|m|/n$ is not close to zero (see the discussion following Theorem 3.1).

The impact of Theorem 4.1 for inference after model selection is immediate in view of Bonferroni's inequality. For the following results, consider the collection $\mathcal{M}$ of candidate models introduced in Section 2. Recall that $\#\mathcal{M}$ and $|\mathcal{M}|$ denote the total number of candidate models and the number of parameters in the most complex candidate model, respectively; moreover, recall that $m_\rho$ denotes the best candidate model and $\hat{m}$ denotes the empirically best candidate model [in the sense that they minimize $\rho^2(m)$ and $\hat{\rho}^2(m)$, respectively, over $m \in \mathcal{M}$].

COROLLARY 4.2. *For $\varepsilon$ satisfying $0 < \varepsilon \le \log(2)$, and uniformly over all data-generating processes as in (2), we have*

$$P\bigg(\|\hat{\mathbb{L}}(\hat{m}) - \mathbb{L}(\hat{m})\|_{\text{TV}} > \frac{1}{\sqrt{n}} + \varepsilon\bigg) \le 7\exp\bigg[\log \#\mathcal{M} - \frac{n-|\mathcal{M}|}{2}\frac{\varepsilon^2}{\varepsilon+2}\bigg].$$

During the proof, we actually derive a slightly stronger version of Corollary 4.2. The result continues to hold if the left-hand side in the preceding inequality is replaced by $P(\max_{m \in \mathcal{M}} \|\hat{\mathbb{L}}(m) - \mathbb{L}(m)\|_{\mathrm{TV}} > 1/\sqrt{n} + \varepsilon)$. This can be used, say, to conduct inference based on the model selected by another model selection procedure that differs from $\hat{m}$.

For the rest of this section, we illustrate the use of our results to construct symmetric prediction intervals centered at $\hat{y}^{(f)}(\hat{m})$ that are approximately valid and short. Similar results can be obtained for one-sided prediction intervals or for testing whether, say, the future response lies above (or below) a prespecified value. Conditional on the training sample, the prediction error $\hat{y}^{(f)}(\hat{m}) - y^{(f)}$ is distributed as $\mathbb{L}(m) \equiv N(\nu(\hat{m}), \delta^2(\hat{m}))$. Hence, a "prediction interval" for $y^{(f)}$ with conditional coverage probability $1 - \alpha$ is given by $[\hat{y}^{(f)}(\hat{m}) - \nu(\hat{m}) - q_\alpha \delta(\hat{m}), \hat{y}^{(f)}(\hat{m}) - \nu(\hat{m}) + q_\alpha \delta(\hat{m})]$, and we write this "prediction interval" informally as

$$\hat{y}^{(f)}(\hat{m}) - \nu(\hat{m}) \pm q_\alpha \delta(\hat{m});$$

here, $q_\alpha$ is the $1 - \alpha/2$ quantile of the standard normal distribution. Note that this construction is infeasible, because it depends on unknown parameters through $\nu(\hat{m})$ and $\delta(\hat{m})$. Corollary 4.2 suggests that a feasible prediction interval can be obtained by replacing the true distribution $\mathbb{L}(\hat{m})$ by the approximating distribution $\hat{\mathbb{L}}(\hat{m})$ and constructing a prediction interval with nominal coverage probability $1 - \alpha$ using $\hat{\mathbb{L}}(\hat{m})$. This amounts to replacing $\nu(\hat{m})$ and $\delta(\hat{m})$ by zero and by $\hat{\delta}(\hat{m})$, respectively, in the preceding display. The resulting prediction interval will be denoted by $\mathcal{I}(\hat{m})$ and is given by

(8) $$\mathcal{I}(\hat{m}) : \hat{y}^{(f)}(\hat{m}) \pm q_\alpha \hat{\delta}(\hat{m}).$$

In view of Corollary 4.2, we get the following result.

PROPOSITION 4.3. *Fix $\varepsilon$ satisfying $0 < \varepsilon \le \log(2)$. Conditional on the training sample, the coverage probability of the prediction interval $\mathcal{I}(\hat{m})$ is within $1/\sqrt{n} + \varepsilon$ of the nominal level, that is,*

$$|(1 - \alpha) - P(y^{(f)} \in \mathcal{I}(\hat{m})|Y, X)| \le \frac{1}{\sqrt{n}} + \varepsilon,$$

*except on an event whose probability is not larger than*

$$7 \exp\left[\log \#\mathcal{M} - \frac{n - |\mathcal{M}|}{2} \frac{\varepsilon^2}{\varepsilon + 2}\right],$$

*uniformly over all data-generating processes as in (2).*

The (infeasible) valid prediction interval based on the selected model $\hat{m}$ discussed prior to Proposition 4.3 has width $2q_\alpha \delta(\hat{m})$, and the width of the feasible interval $\mathcal{I}(\hat{m})$ is $2q_\alpha \hat{\delta}(\hat{m})$. From the perspective of interval width, the "best" model is $m_\delta$, that is, the model minimizing $\delta^2(m)$ over $m \in \mathcal{M}$ (see the discussion at the end of Section 2), and the corresponding exact shortest "prediction interval" is

$$\hat{y}^{(f)}(m_\delta) - \nu(m_\delta) \pm q_\alpha \delta(m_\delta). \tag{9}$$

Again, this construction is infeasible because $m_\delta$ and also $\nu(m_\delta)$ and $\delta(m_\delta)$ depend on unknown parameters. The following result compares the feasible interval $\mathcal{I}(\hat{m})$ based on the selected model with the infeasible shortest possible interval (9) in terms of width, by comparing $\hat{\delta}(\hat{m})$ and $\delta(m_\delta)$.

PROPOSITION 4.4. *For each $\varepsilon > 0$ and uniformly over all data-generating processes as in (2), we have*

$$P\left(\left|\log \frac{\hat{\delta}(\hat{m})}{\delta(m_\delta)}\right| > \varepsilon\right) \leq 4\exp\left[\log \#\mathcal{M} - \frac{n - |\mathcal{M}|}{2}\frac{\varepsilon^2}{\varepsilon + 2}\right].$$

If the upper bound in Proposition 4.4 is small, the length of the feasible prediction interval $\mathcal{I}(\hat{m})$ is close to the length of the (infeasible) shortest prediction interval (9) with high probability. Together with Proposition 4.3, this result establishes that the interval $\mathcal{I}(\hat{m})$ is approximately valid and short with high probability, provided only $n - |\mathcal{M}|$ is large enough compared to $\log \#\mathcal{M}$ (i.e., provided only that the degrees of freedom in the most complex candidate model is large compared to the logarithm of the number of candidate models).

**5. Simulation example.** We now present an example where we search for a "sparse" model in a pool of more than $10^{15}$ candidate models using a training sample of 2000 observations. We demonstrate that a good candidate model can be identified, that the performance of the selected model can be estimated with reasonable accuracy and that a prediction interval post model selection obtains an actual coverage probability reasonably close to its nominal one. The example is meant for demonstration only and should not be mistaken for an exhaustive simulation study (see also [20] for related simulations also covering non-Gaussian scenarios).

Consider a situation where we have available a training sample of 2000 independent observations of the response $y$ and of 1000 explanatory variables $x_j$, $j = 1, \ldots, 1000$, from (2), and where we suspect that the corresponding first 1000 coefficients of $\beta$ are "sparse" in the sense that most of them are very small or zero while a few groups of adjacent coefficients are large. To come up with a collection of candidate models that can pick out the suspected groups of "important" coefficients while not being too large (in the sense that $\log \#M \ll n - |\mathcal{M}|$), we divide the first 1000 coefficients of $\beta$ into 50 blocks of length 20 each, and we consider all candidate models that include or exclude a block at a time (plus the intercept that is always included). This gives $2^{50}$ or a little over $10^{15}$ candidate models.[9] With an exhaustive search over this model space being infeasible, we resort to the obvious greedy general-to-specific strategy. We fit the "overall" model containing all 50 blocks, and eliminate that block whose elimination leads to the smallest increase in the residual sum of squares. This results in a model containing 49 blocks, and now we proceed inductively until all blocks have been eliminated and only the intercept remains.[10] This results in a data-driven rearrangement of the blocks and, thus, of the whole parameter vector $\beta$. The selected model here is the minimizer of $\hat{\rho}^2(\cdot)$ among the models visited by the greedy search and will be denoted by $\hat{m}_g$ throughout this section.

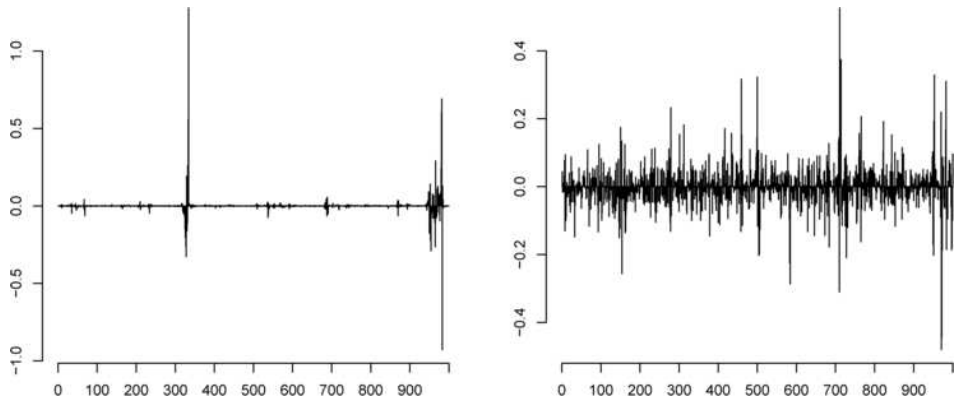

FIG. 1. *The first 1000 coefficients of $\beta$ in the sparse case (left panel) and in the nonsparse case (right panel).*

The suspicion that $\beta$ is sparse, which motivated our choice of candidate models, may or may not be correct in practice. For the true value of the parameter $\beta$, we therefore consider two scenarios. In the first, the coefficients of $\beta$ are indeed sparse, and, in the second, they are not (i.e., a "sparse" and a "nonsparse" case). The first 1000 coefficients of $\beta$ in both the sparse and the nonsparse case are displayed in Figure 1; the remaining coefficients of $\beta$ are set to zero. The first 1000 coefficients of $\beta$ were obtained from realizations of ARCH-processes with different parameters for the sparse and for the non-sparse case.[11] In both cases, the coefficients of $\beta$ were also scaled so that the "signal-to-noise" ratio is five in the sense that $(\text{Var}(y) - \text{Var}(u))/\text{Var}(u) = 5$. If the signal-to-noise ratio is too small, only very parsimonious models perform well; a large signal-to-noise ratio has the opposite effect. We chose a signal-to-noise ratio between these two extremes. The remaining parameters in (2) were chosen as follows. We chose $u$ independent of the $x_j$'s with mean 0 and variance 1, the variance/covariance structure of the explanatory variables was chosen so that $\text{Cov}(x_j, x_k) = 2^{-|j-k|}$ for $j, k \in \{2, \ldots, 1000\}$ (recall that $x_1 = 1$ is the intercept), and we took independent realizations of a standard normal for the means of the $x_j$'s scaled so $\sum_{j=2}^{1000} E(x_j)\beta_j = \sqrt{2}$.[12]

---

[9]We have also experimented with larger (smaller) block-sizes that lead to correspondingly smaller (larger) classes of candidate models. Larger block-sizes give better accuracy of $\hat{\rho}^2(\cdot)$ as an estimator for $\rho^2(\cdot)$ and better coverage properties of prediction intervals post model selection; smaller block-sizes have the opposite effect.

[10]The study of alternative and potentially superior strategies of searching through model space is beyond the scope of this paper.

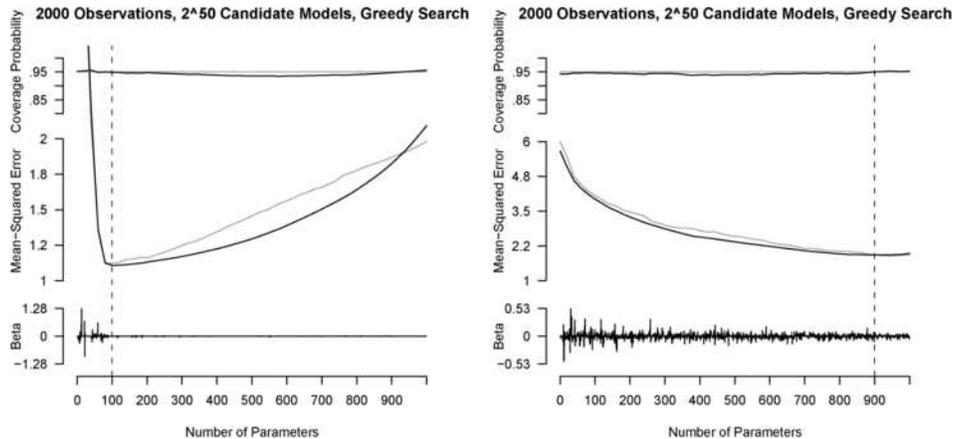

Fig. 2. *Results from one simulation run for the sparse case (left panel) and the nonsparse case (right panel). The graphics are described in the main text.*

For one set of training data [i.e., for 2000 observations from (2) with the parameters just described] the results in both the sparse and the nonsparse case are visualized in Figure 2. The data-driven rearrangement of $\beta$ obtained by the greedy search is shown at the bottom of each panel next to the axis labeled Beta. The block of 20 coefficients to the far right was eliminated first, the block next to it was eliminated next, et cetera, until only the intercept remained. Note that this corresponds to a data-driven sequence of 51 nested models of increasing complexity, from the model containing only the intercept up to the overall model containing all 1000 explanatory variables; the horizontal axis can be thought of as indexing these 51 nested models.

---

[11] In additional experiments, with several other choices for $\beta$, we obtained results consistent with those presented here.

[12] Our choice for the variance/covariance structure of the $x_j$'s is ad-hoc. Repeating the simulations with $\text{Cov}(x_j, x_k) = r^{|j-k|}$ for other values of $r$ between 0 and 0.9, we obtained basically identical results. The same applies, mutatis mutandis, to the choice of the $E(x_j)$'s and to the scaling of the means.

The performance of each of the 51 models obtained by the greedy search is shown by the graph in the middle of each panel, next to the axis labeled Mean-Squared Error. The black line shows the value of $\hat{\rho}^2(\cdot)$ for each of the models (estimated performance), while the gray line shows the true value of $\rho^2(\cdot)$ (actual performance). For better readability, points are joined by lines. The selected model is indicated by a vertical dashed line. Note that the conditional mean-squared error of any predictor is bounded from below by $\text{Var}(u)$, which equals 1 here. Therefore, models $m$ with $\rho^2(m)$ close to 1 perform well. Finally, for each of the 51 models obtained through the greedy search, the coverage probability of a prediction interval based on that model is shown at the top of each panel, next to the axis labeled Coverage Probability. For each such model $m$, we computed the prediction interval $\mathcal{I}(m)$ introduced in (8) with $1-\alpha = 0.95$. The black line shows the true conditional coverage probability of these intervals, conditional on the training sample. Again, points are joined by lines. The gray horizontal line at 0.95 is for reference. Note that models with small $\hat{\rho}^2(m)$ also correspond to short prediction intervals, because the width of $\mathcal{I}(m)$ is governed by $\hat{\delta}(m) \equiv \hat{\rho}(m)$. The conditional coverage probability of the prediction interval based on the selected model is indicated by the vertical dashed line. Because coverage probabilities are computed conditional on the training sample, they can be both above and below the nominal value of 0.95. Because the 51 models shown in each panel of Figure 2 were obtained through a greedy search through model space, $\hat{\rho}^2(\cdot)$ tends to under-estimate $\rho^2(\cdot)$ for these models, resulting in prediction intervals that tend to be too short and whose conditional coverage probabilities tend to fall below 0.95.

In the sparse case (left panel), the chosen class of candidate models is satisfactory in the sense that it contains a relatively parsimonious candidate model that performs well. The selected model $\hat{m}_g$ contains 100 explanatory variables [and coincides with the model minimizing the actual performance $\rho^2(\cdot)$ among the 51 candidates identified by the greedy search]. The selected model's estimated performance of $\hat{\rho}^2(\hat{m}_g) = 1.110$ is close to its actual performance of $\rho^2(\hat{m}_g) = 1.124$ (which in turn is close to the lower bound 1). The conditional coverage probability of the prediction interval based on the selected model is 0.948. In the nonsparse case (right panel), the class of candidate models is unsatisfactory in the sense that it does not contain a simple model that performs well. Among the 51 models identified by the greedy search, the model with 940 explanatory variables performs best [minimizer of $\rho^2(\cdot)$], while the model $\hat{m}_g$ selected by minimizing $\hat{\rho}^2(\cdot)$ contains 900 coefficients. The actual performance of the selected model is $\rho^2(\hat{m}_g) = 1.929$, while its estimated performance is $\hat{\rho}^2(\hat{m}_g) = 1.924$. The selected model improves little over the overall model containing all 1000 explanatory variables, in terms of actual performance as well as in terms of estimated performance.

The conditional coverage probability of the prediction interval based on the selected model is 0.949. Overall, the minimum of the coverage probabilities over the 51 models found by the greedy search (i.e., the minimum of the curve next to the axis labeled `Coverage Probability`) is 0.935 in the sparse case and 0.938 in the nonsparse case, respectively.

The experiment whose results are shown in Figure 2 was repeated a total of 100 times. The results of these repetitions are so similar to those shown in Figure 2 that we do not present them here in detail. Over the 100 repetitions and for the conditional coverage probability corresponding to the selected model, we obtained a median of 0.949 and a minimum of 0.936 in the sparse case, and a median of 0.942 and a minimum of 0.924 in the nonsparse case. Also, over 100 repetitions and for the minimum of the coverage probabilities corresponding to the 51 models identified by the greedy search, we obtained a median of 0.931 and a minimum of 0.919 in the sparse case and a median of 0.934 and a minimum of also 0.918 in the nonsparse case.

## 6. Remarks and extensions.

REMARK 6.1 [*Examples of data-generating processes as in* (2)]. The distribution of the random variables in (2) is, of course, characterized by their first and second moments. Assume, for simplicity, that the $x_j$'s are uncorrelated with $u$ and that $u$ has mean zero. Write $\beta$ for the sequence of regression coefficients $\beta = (\beta_j)_{j \geq 1}$, write $\sigma^2$ for the variance of $u$ and denote the sequence of means and the variance/covariance net of the $x_{j+1}$'s, $j \geq 1$, by $\gamma = (\gamma_j)_{j \geq 1}$ and $\Sigma = (\Sigma_{j,k})_{j,k=1}^{\infty}$ (recall that $x_1$ denotes the intercept, i.e., $x_1 = 1$). That is, the mean and the variance of $x_{j+1}$ are $\gamma_j$ and $\Sigma_{j,j}$, respectively, and the covariance of $x_{j+1}$ and $x_{k+1}$ is $\Sigma_{j,k}$, $1 \leq j < k$. Then, the (joint) distribution of $y$, $x_j$, $j \geq 1$ and $u$ in (2) is characterized by $(\beta, \gamma, \Sigma, \sigma)$. Write $\Xi$ for the collection of all quadruples $\xi = (\beta, \gamma, \Sigma, \sigma)$ such that the series in (2) converges in $L_2$, and such that the joint distribution of the $x_j$'s for $j > 1$ and of $u$ is nondegenerate. The following examples illustrate that $\Xi$ is quite rich and includes subsets that are noncompact (with respect to the appropriate canonical topology):

(i) Assume that the $x_j$'s for $j > 1$ are uncorrelated with common variance equal to unity, and write $I$ for the corresponding variance/covariance net (i.e., $I_{j,k}$ equals one if $j = k$ and zero otherwise). Then, $\Xi$ contains all quadruples $\xi$ of the form $\xi = (\beta, \gamma, I, \sigma)$ satisfying $\beta \in l_2$, $\gamma \in l_2$, and $\sigma > 0$.

(ii) Let $\varsigma = (\varsigma_j)_{j \geq 1}$ be an arbitrary sequence of positive numbers. Assume that the $x_j$'s are uncorrelated as in (i) before but now with $\text{Var}(x_j) = \varsigma_j^2$, and write $\text{diag}(\varsigma^2)$ for the corresponding variance/covariance net. Then, $\Xi$ contains all quadruples $\xi$ of the form $\xi = (\beta, \gamma, \text{diag}(\varsigma^2), \sigma)$ for which $\beta\gamma \in l_1$ and $\beta\varsigma \in l_2$ (where the products are understood component-wise) and $\sigma > 0$.

(iii) Fix $p$ satisfying $1 \leq p \leq \infty$, and let $q$ be such that, either $1 < p < \infty$ and $1/p + 1/q = 1$, or $p = 1$ and $q = \infty$, or $p = \infty$ and $q = 1$. Let $S : l_p \to l_q$ be a continuous linear operator satisfying $\langle \alpha, S\beta \rangle = \langle S\alpha, \beta \rangle$ for each $\alpha$ and $\beta$ in $l_p$, and satisfying $\langle \alpha, S\alpha \rangle > 0$ whenever $\alpha \in l_p$ is nonzero. Here, $\langle \cdot, \cdot \rangle$ denotes the usual product of sequences (i.e., the sum of component-wise products). The operator $S$ defines a variance/covariance net $\Sigma(S)$ by $\Sigma(S)_{j,k} = \langle e_j, Se_k \rangle$, where $e_l$ denotes a sequence with a 1 in the $l$th position and zeroes otherwise ($l \geq 1$). Then, $\Xi$ contains all quadruples $\xi$ of the form $\xi = (\beta, \gamma, \Sigma(S), \sigma)$ satisfying $\beta \in l_p$, $\gamma \in l_q$, $S$ as before, and $\sigma > 0$.

REMARK 6.2 (*Reduced-rank models*). We have required that the joint distribution of the $x_j$'s for $j > 1$ in (2) is nondegenerate (and Gaussian). For our purpose, this guarantees that the $n \times |m|$ matrix of those regressors in the training sample that are included in a model $m \in \mathcal{M}$ is nondegenerate with probability one. We now discuss the case where this requirement is not met. Assume, for a candidate model $m \in \mathcal{M}$, that some of the explanatory variables $x_j$ that are included in the model $m$ are perfectly correlated with each other. In that case, there is a submodel $m'$ of $m$ (i.e., a model $m'$ satisfying $m'_j \leq m_j$ for each $j$), such that the explanatory variables included in model $m'$ are not perfectly correlated with each other, and such that the least-squares predictors based on model $m$ and $m'$ coincide [i.e., $\hat{y}^{(f)}(m) = \hat{y}^{(f)}(m')$], almost surely. Here, the restricted least-squares estimator $\hat{\beta}(m)$ needs to be computed using a generalized inverse in the least-squares formula, because the sample regressor matrix corresponding to model $m$ is of reduced rank, almost surely. Hence, we also have $\rho^2(m) = \rho^2(m')$ and $\mathbb{L}(m) \equiv \mathbb{L}(m')$ a.s. Now, repeat this replacement process for each candidate model in $\mathcal{M}$ (i.e., replace each reduced-rank model by an appropriate full-rank submodel and leave the full-rank models unchanged). This results in a new collection of candidate models, which we denote by $\mathcal{M}'$. Inspection of the proofs reveals that all the results in Sections 2, 3 and 4 continue to hold with $\mathcal{M}'$ replacing $\mathcal{M}$.

REMARK 6.3 (*Note on constants*). Several performance bounds that are reported in this paper depend on the constants $\#\mathcal{M}$ and $|\mathcal{M}|$ (see Corollaries 3.2 and 4.2, as well as Propositions 4.3 and 4.4). These bounds are conservative because they hold uniformly over a large class of data-generating processes and for each class $\mathcal{M}$ of candidate models that satisfies (3). In particular, the results also cover the case where all candidate models are equally complex and where all fit equally well. In view of this, it is not surprising to find the constants $\#\mathcal{M}$ and $|\mathcal{M}|$ in the upper bounds. If additional regularity conditions are imposed on the regression parameter, and if the family of candidate models $\mathcal{M}$ is chosen in accordance to these regularity conditions (e.g., sparse candidate models in case a sparsity condition is imposed on the true regression parameter), it is likely that the upper bounds can be improved. Also, the fact that the upper bounds all increase linearly with the number of candidate models (i.e., with $\#\mathcal{M}$) originates in the use of Bonferroni's inequality, which could leave room for improvement. These issues, however, are beyond the scope of this paper.

REMARK 6.4 (*Asymptotic rates*). The results in Sections 3 and 4 allow us to read off the rates at which quantities like $\log \hat{\rho}^2(m)/\rho^2(m)$ or $\log \hat{\rho}^2(\hat{m})/\rho^2(m_\rho)$, for example, converge to zero, in probability, in appropriate asymptotic settings. Under rather weak conditions, we show that the typical rate is $1/\sqrt{n}$ in the following:

(i) Consider a sequence of sample sizes $n$, and a corresponding sequence of candidate models $m^{(n)}$ (that may depend on $n$), such that $|m^{(n)}|/n \leq r$ for fixed $r$, $0 < r < 1$, and for each $n$. As always, we also assume that $m_1^{(n)} = 1$ and that $|m^{(n)}| < n - 1$. Denoting the distribution of the sample of size $n$ by $P_n(\cdot)$, Theorem 3.1 entails that

$$P_n\bigg(\sqrt{n}\bigg|\log\frac{\hat{\rho}^2(m^{(n)})}{\rho^2(m^{(n)})}\bigg| > t\bigg) \leq 6\exp\bigg[-\frac{1-r}{8}\frac{t^2}{t+8}\bigg],$$

for each $n$, for each $t > 0$, and uniformly over all data-generating processes as in (2). In other words, $\log \hat{\rho}^2(m^{(n)})/\rho^2(m^{(n)})$ is of order $1/\sqrt{n}$ in probability, uniformly over all data-generating processes as in (2). In a similar fashion, $\|\hat{\mathbb{L}}(m^{(n)}) - \mathbb{L}(m^{(n)})\|_{\mathrm{TV}}$ is uniformly of order $1/\sqrt{n}$ in probability (see Theorem 4.1).

(ii) Now, consider a sequence of sample sizes $n$ and a corresponding sequence of families of candidate models $\mathcal{M}^{(n)}$ [such that (3) holds for each $m \in \mathcal{M}^{(n)}$ and for each $n$]. Moreover, assume that $|\mathcal{M}^{(n)}| < r$ for fixed $r$, $0 < r < 1$, and for each $n$, and that $\log \#\mathcal{M}^{(n)} = o(n)$. We stress that now quantities like the "best" model $m_\rho$, the empirically best model $\hat{m}$, the conditional distribution of the prediction error $\mathbb{L}(m)$, its estimated version $\hat{\mathbb{L}}(m)$, the prediction interval $\mathcal{I}(\hat{m})$, et cetera, all depend on $n$, although this dependence is not shown explicitly by the notation. Under these assumptions we obtain that the following quantities are each of order $1/\sqrt{n}$ in probability, uniformly over all data-generating processes as in (2): $\log \rho^2(\hat{m})/\rho^2(m_\rho)$ and $\log \hat{\rho}^2(\hat{m})/\rho^2(\hat{m})$ (see Corollary 3.2); $\|\hat{\mathbb{L}}(\hat{m}) - \mathbb{L}(\hat{m})\|_{\mathrm{TV}}$ (see Corollary 4.2); $(1-\alpha) - P_n(y^{(f)} \in \mathcal{I}(\hat{m}_n)|Y,X)$ (see Proposition 4.3); and $\log \hat{\delta}^2(\hat{m})/\delta^2(m_\delta)$ (see Proposition 4.4).

## APPENDIX A: PROOFS FOR SECTION 2

The following two lemmas will be instrumental in the proof of Proposition 2.1. We suspect that these two results, which basically rely on the rotational invariance of the normal distribution, are well known, in some form or another, but we could not find a convenient reference in the literature. Throughout, the Euclidean norm of a vector $v \in \mathbb{R}^k$ is denoted by $\|v\|$.

LEMMA A.1. *For $k \geq 1$, let $a \sim N(0, I_k)$ and fix $a_0 \in \mathbb{R}^k$ with $\|a_0\| = 1$. Then, there exists a $k \times k$ matrix $R$ whose elements are measurable functions of $a$, such that $R'R = I_k$ and $a = \|a\| R a_0$ almost surely.*

PROOF. Write $e_j$ for the $j$th Euclidean basis vector of $\mathbb{R}^k$. It suffices to consider the case where $a_0 = e_1$, because $a_0$ can be written as $a_0 = S e_1$ where $S$ is a fixed orthonormal $k \times k$ matrix. For $a_0 = e_1$, consider the event $E$ where $\|a\| > 0$ and where $a$ is linearly independent of $e_2, \ldots, e_k$. Clearly, $E$ is an almost sure event. On $E$, compute an orthonormal basis $r_1, \ldots, r_k$ of $\mathbb{R}^k$, by setting $r_1 = a/\|a\|$ and by then applying the Gram–Schmidt orthonormalization procedure to $r_1, e_2, \ldots, e_k$, and set $R = (r_1, \ldots, r_k)$. On $E^c$, set $R = I_k$, say. The matrix $R$ has the desired properties. $\square$

LEMMA A.2. *Let $M$ be a $k \times l$ matrix with i.i.d. standard normal entries, and let $a_0 \in \mathbb{R}^k$ with $\|a_0\| = 1$ ($1 \leq l \leq k$). For $P_M = M(M'M)^{-1}M'$, the distribution of $a_0' P_M a_0$ is given by*

$$a_0' P_M a_0 \sim \frac{\chi_l^2}{\chi_l^2 + \chi_{k-l}^2},$$

*where $\chi_l^2$ and $\chi_{k-l}^2$ denote two independent chi-square random variables with the indicated degrees of freedom. In case $l = k$, $\chi_{k-l}^2$ is to be interpreted as constant equal to zero and the distribution on the right-hand side of the preceding display is to be interpreted as point mass at one.*

PROOF. As the case $l = k$ is trivial, we may assume that $l < k$. Let $a \sim N(0, I_k)$ independent of $M$. Using Lemma A.1, we can rewrite $a$ as $a = \|a\| R a_0$ almost surely. With probability one, we thus have $a_0 = R'a/\|a\|$ and $a_0' P_M a_0$ can be written as

$$a_0' P_M a_0 = \frac{a' R P_M R' a}{a'a} = \frac{a' RM(M'R'RM)^{-1}M'R'a}{a'a} = \frac{a' P_{M_\circ} a}{a'a}$$

almost surely, where, for the last equality, we use $M_\circ$ as shorthand for $RM$ and define $P_{M_\circ}$ like $P_M$ with $M_\circ$ replacing $M$. Conditional on $a$, the columns of $M_\circ$ are i.i.d. $N(0, I_k)$ [because the columns of $M$ are i.i.d. $N(0, I_k)$ independent of $a$, and because $RR' = I_k$]. As that conditional distribution does not depend on the conditioning variable, we see that the entries of $M_\circ$ are i.i.d. standard Gaussians independent of $a$. Hence, the columns of $M_\circ$ are linearly independent with probability one, and the expression on the far right-hand side of the preceding display is distributed as $\chi_l^2/(\chi_l^2 + \chi_{k-l}^2)$.
□

Before turning to the proof of Proposition 2.1, the following preparatory consideration and the attending lemma are also required. Throughout the following, fix a candidate model $m \in \mathcal{M}$. Recall the linear model (2), and write $z$ for the $|m|$-vector of those explanatory variables $x_j$ that are included in the model $m$ (in their natural order, so that $z_1$ corresponds to the intercept, i.e., $z_1 = x_1 = 1$). Because $y$ and $z$ are jointly Gaussian, the conditional distribution of $y$ given $z$ is again a Gaussian. Because the model $m$ includes an intercept (i.e., $z_1 = 1$) the conditional mean of $y$ given $z$ is a linear function of $z$. Recalling that the conditional variance of $y$ given $z$ is $\sigma^2(m)$, we see that $y|z \sim N(z'\theta, \sigma^2(m))$ for an appropriate $|m|$-vector $\theta$. In other words, (2) can be rewritten as

$$(10) \qquad y = z'\theta + v$$

with $v \sim N(0, \sigma^2(m))$ independent of $z$. The vector $z$ of those explanatory variables that are included in model $m$ is also Gaussian, and, in the following, we write $\eta$ and $\Gamma$ for the mean-vector and for the variance/covariance matrix of its distribution, respectively:

$$(11) \qquad z \sim N(\eta, \Gamma).$$

Clearly, $\eta$ is an $|m|$-vector and $\Gamma$ is an $|m| \times |m|$ matrix. Because the first regressor corresponds to the intercept, we have $\eta_1 = 1$ and $\Gamma_{1,1} = 0$. In case $|m| > 1$, the submatrix of $\Gamma$ corresponding to $z_2, \ldots, z_{|m|}$ is positive definite by assumption [see the discussion following (2)].

LEMMA A.3. *For fixed $m \in \mathcal{M}$, let $\eta$ and $\Gamma$ be as in (11) and set $\Delta = \Gamma + \eta\eta'$. Then, $\Delta$ is positive definite, and so is its symmetric square root $\Delta^{1/2}$. Moreover, the matrix $\Delta^{-1/2}\Gamma\Delta^{-1/2}$ admits a spectral representation $\Delta^{-1/2}\Gamma\Delta^{-1/2} = W\Lambda W'$ such that $\Lambda = \mathrm{diag}(0, 1, \ldots, 1)$ (i.e., the first eigenvalue equals zero and all the others equal one), such that $W = (w^{(1)}, \ldots, w^{(|m|)})$ with the $w^{(j)}$, $j = 1, \ldots, |m|$, being orthonormal eigenvectors, and such that $w^{(1)} = \Delta^{-1/2}\eta$. In particular, $W'\Delta^{-1/2}\eta = (1, 0, \ldots, 0)' \in \mathbb{R}^{|m|}$.*

PROOF. To show that $\Delta > 0$, assume that $w \in \mathbf{R}^{|m|}$ is such that $w'\Delta w = 0$. Partition $w$ as $w = (w_1, w'_{\neg 1})'$ (i.e., into its first component $w_1$ and the vector $w_{\neg 1}$ of its $|m| - 1$ remaining components), partition $\eta$ conformably as $\eta = (\eta_1, \eta'_{\neg 1})'$, and let $\Sigma$ denote the lower diagonal $(|m|-1) \times (|m|-1)$ submatrix of $\Gamma$. Recall that $\eta_1 = 1$, that $\Sigma > 0$, and that the first row as well as the first column of $\Gamma$ contain zeroes only [see (11) and the attending discussion]. Therefore,

$$w'\Delta w = w'\Gamma w + w'\eta\eta' w = w'_{\neg 1}\Sigma w_{\neg 1} + w'\eta\eta' w.$$

Because $w'\Delta w = 0$ and $\Sigma > 0$, we see that $w_{\neg 1} = 0$. Hence, $w'\Delta w = w_1^2\eta_1^2 = w_1^2$, so that $w_1$ also equals zero and $w = 0$.

Write $K$ as shorthand for $\Delta^{-1/2}\Gamma\Delta^{-1/2}$, and note that $K$ has rank $|m|-1$, because $\Gamma$ has rank $|m|-1$. Moreover, we have

$$K = \Delta^{-1/2}\Gamma\Delta^{-1/2} = \Delta^{-1/2}(\Delta - \eta\eta')\Delta^{-1/2} = I_{|m|} - \Delta^{-1/2}\eta\eta'\Delta^{-1/2}.$$

Set $w^{(1)} = \Delta^{-1/2}\eta$, and note that $w^{(1)}$ is nonzero because $\eta_1 = 1$. For each vector $w$ in the orthogonal complement of $w^{(1)}$, we thus have $Kw = w$. Hence, $|m|-1$ eigenvalues of $K$ equal one, and the corresponding eigenvectors (which can be chosen as normalized and mutually orthogonal) are orthogonal to $w^{(1)}$. The remaining eigenvalue of $K$ must be zero, and $w^{(1)}$ must be a corresponding eigenvector. This entails that $0 = Kw^{(1)} = w^{(1)} - w^{(1)}w^{(1)'}w^{(1)}$, whence $\|w^{(1)}\| = 1$. Finally, $W'\Delta^{-1/2}\eta = W'w^{(1)} = (1, 0, \ldots, 0)'$. $\square$

PROOF OF PROPOSITION 2.1. Without loss of generality, we may assume that the random matrices in the following arguments are invertible whenever we need them to be, because the event where that is not the case has probability zero. For the given model $m$, write $Z$ for the $n \times |m|$ matrix of those explanatory variables in the training sample that are included in the model $m$, such that the $i$th entry of $Y$ and the $i$th column of $Z'$ are independent copies of $y$ and $z$ in (10) for $i = 1, \ldots, n$. Note that the $i$th column of $Z'$ is distributed as in (11). Let $\Delta^{1/2}$, $W$ and $\Lambda$ be as in Lemma A.3, and set $Z^{(\bullet)} = Z\Delta^{-1/2}W$. Lemma A.3 now entails that the $i$th column of $Z^{(\bullet)'}$ is distributed as $N(e_1, \Lambda)$, where $e_1 = (1, 0, \ldots, 0)' \in \mathbb{R}^{|m|}$ and $\Lambda$ is the diagonal matrix $\Lambda = \text{diag}(0, 1, \ldots, 1)$. In particular, $Z^{(\bullet)}$ can be partitioned as $Z^{(\bullet)} = (\iota, Z^{(\circ)})$, where $\iota$ is an $n$-vector of ones and $Z^{(\circ)}$ is an $n \times (|m|-1)$ matrix with i.i.d. standard normal entries.

For $\theta$ as in (10), set $V = Y - Z\theta$, and note that $V \sim N(0, \sigma^2(m)I_n)$ independent of $Z$. From this, it follows that $\hat{\sigma}^2(m) \sim \sigma^2(m)\chi^2_{n-|m|}/(n-|m|)$ as claimed. Moreover, $\nu(m)$ and $\delta^2(m)$ can be written as

$$\nu(m) = \eta'(Z'Z)^{-1}Z'V \quad \text{and} \quad \delta^2(m) = V'Z(Z'Z)^{-1}\Gamma(Z'Z)^{-1}Z'V + \sigma^2(m)$$

[compare the definitions of $\nu(m)$ and $\delta^2(m)$ in Section 2, as well as (10) and (11)]. From the first equation in the preceding display, we also see that $E[\nu(m)] = 0$. It remains to derive the distribution of $\nu^2(m)$ and of $\delta^2(m)$.

To this end, we need more convenient representations of these quantities. Rewrite $\nu(m)$ as

$$\begin{aligned}\nu(m) &= \eta'\Delta^{-1/2}(\Delta^{-1/2}Z'Z\Delta^{-1/2})^{-1}\Delta^{-1/2}Z'V \\ &= \eta'\Delta^{-1/2}W(W'\Delta^{-1/2}Z'Z\Delta^{-1/2}W)^{-1}W'\Delta^{-1/2}Z'V \\ &= e_1'(Z^{(\bullet)'}Z^{(\bullet)})^{-1}Z^{(\bullet)'}V,\end{aligned}$$

where the last equality follows upon observing that we have set $Z^{(\bullet)} = Z\Delta^{-1/2}W$ and that $W'\Delta^{-1/2}\eta = e_1$ by Lemma A.3. A similar argument gives

$$\begin{aligned}\delta^2(m) - \sigma^2(m) &= V'Z^{(\bullet)}(Z^{(\bullet)'}Z^{(\bullet)})^{-1}W'\Delta^{-1/2}\Gamma\Delta^{-1/2}W(Z^{(\bullet)'}Z^{(\bullet)})^{-1}Z^{(\bullet)'}V \\ &= V'Z^{(\bullet)}(Z^{(\bullet)'}Z^{(\bullet)})^{-1}\Lambda(Z^{(\bullet)'}Z^{(\bullet)})^{-1}Z^{(\bullet)'}V,\end{aligned}$$

where we use the spectral representation of $\Delta^{-1/2}\Gamma\Delta^{-1/2}$ given in Lemma A.3 to get the last equality. We thus see that $\nu^2(m)$ is the square of the first component of the $|m|$-vector $(Z^{(\bullet)'}Z^{(\bullet)})^{-1}Z^{(\bullet)'}V$, and $\delta^2(m) - \sigma^2(m)$ is the sum of squares of the remaining $|m|-1$ components of that vector.

Partitioning $Z^{(\bullet)}$ as $Z^{(\bullet)} = (\iota, Z^{(\circ)})$ as before, we see that

$$(12) \quad (Z^{(\bullet)'}Z^{(\bullet)})^{-1}Z^{(\bullet)'}V = \begin{pmatrix} (\iota'(I_n - P_{Z^{(\circ)}})\iota)^{-1}\iota'(I_n - P_{Z^{(\circ)}})V \\ (Z^{(\circ)'}(I_n - P_\iota)Z^{(\circ)})^{-1}Z^{(\circ)'}(I_n - P_\iota)V \end{pmatrix},$$

where $P_\iota$ and $P_{Z^{(\circ)}}$ denote the orthogonal projections on the space spanned by $\iota$ and on the column space of $Z^{(\circ)}$, respectively. Relation (12) follows either by using the inversion formula for partitioned matrices on the corresponding partition of $Z^{(\bullet)'}Z^{(\bullet)}$ and simplifying, or from geometric properties of orthogonal projections.

For the distribution of $\nu^2(m)$, recall that $\nu^2(m)$ is the square of the first component of the vector on the right-hand side of (12). In particular, $\nu^2(m)$ can be written as

$$\nu^2(m) = \frac{V'P_{(I_n - P_{Z^{(\circ)}})\iota}V}{\iota'(I_n - P_{Z^{(\circ)}})\iota}.$$

The numerator in the preceding display is distributed as $\sigma^2(m)\chi_1^2$, independent of $Z^{(\circ)}$. The denominator is a function of $Z^{(\circ)}$ and, hence, independent of the numerator. Using the Lemma A.2 with $\iota/\sqrt{n}$ and $Z^{(\circ)}$ replacing $a_0$ and $M$ we see, in the notation used in that lemma, that $\iota'P_{Z^{(\circ)}}\iota$ has the same distribution as $n\chi_{|m|-1}^2/(\chi_{|m|-1}^2 + \chi_{n-|m|+1}^2)$. Hence,

$$\iota'(I_n - P_{Z^{(\circ)}})\iota \sim n\frac{\chi_{n-|m|+1}^2}{\chi_{|m|-1}^2 + \chi_{n-|m|+1}^2}.$$

This entails that $\nu^2(m) \sim (\chi_1^2/n)\sigma^2(m)(1+\chi_{|m|-1}^2/\chi_{n-|m|+1}^2)$ as claimed.

For the distribution of $\delta^2(m)$, write $M$ as shorthand for $(I_n - P_\iota)Z^{(\circ)}$. We see from (12) that
$$\delta^2(m) - \sigma^2(m) = V'M(M'M)^{-2}M'V = w'(M'M)^{-1}w,$$
where, for the last equality, we use $w$ to denote the $(|m|-1)$-vector $w = (M'M)^{-1/2}M'V$. Since $V \sim N(0, \sigma^2(m)I_n)$, it follows that $w \sim N(0, \sigma^2(m) \times I_{|m|-1})$, independent of $M$. Using Lemma A.1 with $w$ and $e_1$ replacing $a$ and $a_0$, we obtain an orthonormal matrix $R$ such that $w = \|w\|Re_1$ almost surely. It follows that $\delta^2(m) - \sigma^2(m)$ can be written as
$$\|w\|^2 e_1'(R'MM'R)^{-1}e_1 = \|w\|^2 e_1'(R'Z^{(\circ)'}(I_n - P_\iota)Z^{(\circ)}R)^{-1}e_1.$$
Write $Z^{(R)}$ as shorthand for $Z^{(\circ)}R$. Since $R'R = I_{|m|-1}$, we see that $Z^{(R)} = Z^{(\circ)}R$ is an $n \times (|m|-1)$ matrix with i.i.d. Gaussian entries, and that $Z^{(R)}$ is independent of $w$. Partition $Z^{(R)}$ as $Z^{(R)} = (Z_1^{(R)}, Z_{\neg 1}^{(R)})$, where $Z_1^{(R)}$ is the first column of $Z^{(R)}$, and apply the partitioned inversion formula to the corresponding partition of $(Z^{(R)'}(I_n - P_\iota)Z^{(R)})$. This gives
$$\delta^2(m) - \sigma^2(m) = \frac{\|w\|^2}{Z_1^{(R)'}(I_n - P_\iota)(I_n - P_{(I_n-P_\iota)Z_{\neg 1}^{(R)}})(I_n - P_\iota)Z_1^{(R)}}$$
almost surely. In the preceding expression, the numerator is distributed as $\sigma^2(m)\chi_{|m|-1}^2$ and is independent of the denominator. For the denominator, note that $Z_1^{(R)}$ is an $n$-vector of i.i.d. standard Gaussians, and $(I_n - P_\iota)(I_n - P_{(I_n-P_\iota)Z_{\neg 1}^{(R)}})(I_n - P_\iota)$ is the matrix of an orthogonal projection of rank $n - |m| + 1$ (except on a probability zero event, as is easy to see). It follows that $\delta^2(m) - \sigma^2(m)$ is distributed as $\sigma^2(m)\chi_{|m|-1}^2/\chi_{n-|m|+1}^2$ and $\delta^2(m)$ is distributed as $\sigma^2(m)(1 + \chi_{|m|-1}^2/\chi_{n-|m|+1}^2)$ as required. $\square$

## APPENDIX B: AUXILIARY LEMMAS FOR SECTIONS 3 AND 4

In this section, we show, in essence, that $\delta^2(m)$, $\hat\rho^2(m)$ and $\rho^2(m)$ each are close to the same value with high probability, provided that $n - |m|$ is large enough (see Lemmas B.3, B.4 and B.5, resp.). To derive these results, we also need the two elementary lemmas that follow: Lemma B.1 gives bounds on certain probabilities involving a $\chi_1^2$ random variable, and Lemma B.2 gives a collection of inequalities that will be used later.

LEMMA B.1. *Let $F(\cdot)$ denote the cumulative distribution function (c.d.f.) of the $\chi_1^2$ distribution. Then*
$$F\left(t\frac{\log(t)}{t-1}\right) - F\left(\frac{\log(t)}{t-1}\right) < \frac{\log(t)}{\sqrt{2\pi e}}$$

*holds for each $t > 1$. Moreover, we have*

$$1 - F(t) \leq \sqrt{\frac{2}{\pi}} \exp\left[-\frac{t + \log(t)}{2}\right]$$

*for each $t > 0$.*

PROOF. For the first inequality, write $g(t)$ as shorthand for the left-hand side, and write $h(t)$ for the right-hand side. We need to show that $g(t) \leq h(t)$. Since $\lim_{t \to 1} g(t) = \lim_{t \to 1} h(t) = 0$, this will follow if we can show that $g'(t) \leq h'(t)$. First, note that $g'(t)$ is given by

$$F'\left(t\frac{\log(t)}{t-1}\right)\left[\frac{1}{t-1} - \frac{\log(t)}{(t-1)^2}\right] - F'\left(\frac{\log(t)}{t-1}\right)\left[\frac{1}{t(t-1)} - \frac{\log(t)}{(t-1)^2}\right]$$

$$= \frac{\log(t)}{(t-1)^2}\left[F'\left(\frac{\log(t)}{t-1}\right) - F'\left(t\frac{\log(t)}{t-1}\right)\right] = \frac{1}{\sqrt{2\pi}}\sqrt{\frac{\log(t)}{t-1}}t^{-(2t-1)/(2t-2)},$$

where the two equalities in the preceding display follow by plugging-in the formula $F'(t) = t^{-1/2}\exp(-t/2)/\sqrt{2\pi}$ and simplifying. We need to show that $g'(t)/h'(t) \leq 1$; that is,

$$\sqrt{\frac{\log(t)}{t-1}}e^{-(1/2)(\log(t)/(t-1)-1)} \leq 1$$

[the left-hand side of the preceding inequality equals $g'(t)/h'(t)$, which is easily seen by using the formula for $g'(t)$ obtained before and $h'(t) = 1/(t\sqrt{2\pi e})$]. For $s$ satisfying $0 < s < 1$ set $u(s) = \sqrt{s}\exp(-(s-1)/2)$, and set $v(t) = \log(t)/(t-1)$ for $t$ as before. For each $t > 1$, we have $0 < v(t) < 1$, so that $u(v(t))$ is well defined. Clearly, the left-hand side of the inequality in the preceding display can be written as $u(v(t))$, and we need to show that $u(v(t)) \leq 1$. Since $\lim_{t \to 1} v(t) = 1$ (as is easy to see), we get $\lim_{t \to 1} u(v(t)) = 1$. It hence suffices to show that $u(v(t))$ is decreasing or, equivalently, $\partial u(v(t))/\partial t = u'(v(t))v'(t) \leq 0$ for $t > 1$. It is now elementary to verify that $v'(t) \leq 0$ for $t > 1$ and that $u'(s) > 0$ for $s$ satisfying $0 < s < 1$. Hence, $u'(v(t))v'(t) \leq 0$ and $u(v(t))$ is decreasing.

For the second inequality, write $\Phi(\cdot)$ and $\phi(\cdot)$, respectively, for the c.d.f. and for the Lebesgue density of the standard normal distribution. The result follows upon observing that $1 - F(t) = 2(1 - \Phi(\sqrt{t}))$ and that

$$2(1 - \Phi(\sqrt{t})) \leq 2\frac{\phi(\sqrt{t})}{\sqrt{t}} = \sqrt{\frac{2}{\pi}}\exp\left[-\frac{t + \log(t)}{2}\right],$$

where the inequality holds because of the well-known argument that $1 - \Phi(t) = \int_t^\infty \phi(u)\,du < \int_t^\infty (1 + 1/u^2)\phi(u)\,du = \phi(t)/t$ for $t > 0$. □

LEMMA B.2. (i) *For $s$ satisfying $0 < s < 1$ and for $t \geq 0$, we have*

$$t - s \log \frac{e^t + s - 1}{s} \geq (1-s)\frac{t^2}{t+1+s}.$$

(ii) *For $s$ and $t$ satisfying $0 < s < 1$ and $0 \leq t < -\log(1-s)$, we have*

$$-t - s\log(e^{-t} + s - 1) \geq t - s\log(e^t + s - 1).$$

(iii) *For $t \geq 0$, we have*

$$e^t - 1 - t \geq e^{-t} - 1 + t \geq \frac{t^2}{t+2}.$$

PROOF. For part (i), set $f(t) = t - s\log((e^t + s - 1)/s)$ and $g(t) = (1-s)t^2/(t+1+s)$. To show that $f(t) \geq g(t)$, first note that $f(0)$ and $g(0)$ are both equal to zero. It thus suffices to show that $f'(t) \geq g'(t)$ for each $t > 0$. It is easy to see that

$$f'(t) = (1-s)\frac{e^t - 1}{s + e^t - 1} \quad \text{and} \quad g'(t) = (1-s)\frac{t^2 + 2t(s+1)}{(t+s+1)^2}.$$

Plugging these formulae into the relation $f'(t) \geq g'(t)$ and simplifying, we see that the relation is equivalent to

$$(e^t - 1)(1+s)^2 - st^2 - 2s(1+s)t \geq 0.$$

Replacing $e^t - 1$ by $t + t^2/2$ in the preceding expression, we obtain a lower bound for the left-hand side. After trivial simplifications, that bound reduces to $t(1-s^2) + t^2(1+s^2)/2$ which is nonnegative because $t \geq 0$ and $s \leq 1$.

For part (ii), let $f(t) = t - s\log(e^t + s - 1)$ and $h(t) = f(-t) - f(t)$. We need to show that $h(t) \geq 0$. Since $h(0) = 0$, it remains to show that $h'(t) \geq 0$. Now, $h(t) = -2t - s\log(e^{-t} + s - 1) + s\log(e^t + s - 1)$ and

$$h'(t) = -2 + \frac{se^{-t}}{e^{-t} + s - 1} + \frac{se^t}{e^t + s - 1}.$$

Note that, by choice of $t < -\log(1-s)$ and $t > 0$, both denominators in the two fractions in the preceding display are positive. Multiplying the expressions in the preceding display by $(e^{-t} + s - 1)(e^t + s - 1) > 0$ and simplifying, we see that $h'(t) \geq 0$ if

$$(1-s)(2-s)(e^t - 1)(1 - e^{-t}) \geq 0.$$

This inequality is, of course, always satisfied because $s < 1$ and $t \geq 0$.

For part (iii), first expand $e^t - e^{-t}$ as $\sum_{j=0}^{\infty}(t^j - (-t)^j)/j! = 2\sum_{j=0}^{\infty} t^{2j+1}/(2j+1)!$. Hence, $e^t - e^{-t} \geq 2t$, which is equivalent to the first inequality in (iii). For the second inequality, set $f(t) = e^{-t} - 1 + t$ and $g(t) = t^2/(t+$ 2). Since $f(0) = g(0) = 0$, it suffices to show that $f'(t) \geq g'(t)$, and that inequality is easily seen to be equivalent to

$$4 - e^{-t}(t+2)^2 \geq 0.$$

Write $h(t)$ for the left-hand side of the preceding inequality. Observing that $h(0) = 0$ and that $h'(t) = e^{-t}t(t+2) \geq 0$ completes the proof. $\square$

We are now ready to give the three results that state that $\delta^2(m)$, $\hat{\rho}^2(m)$ and $\rho^2(m)$ each are close to the same fixed value with high probability, provided that $n - |m|$ is large enough. That value is taken as $n\sigma^2(m)/(n-|m|+1)$, which is close but not equal to $E[\delta^2(m)]$ or $E[\rho^2(m)]$ (see the discussion and formula for $E[\rho^2(m)]$ given at the end of Section 2). Throughout, let $m$ be a fixed candidate model from $\mathcal{M}$.

LEMMA B.3.  *For each $t \geq 0$, we have*

$$P\left(\delta^2(m)\frac{n-|m|+1}{n\sigma^2(m)} > \exp(t)\right) \leq \exp\left[-\frac{n-|m|+1}{2}\frac{t^2}{t+1+(|m|-1)/n}\right],$$

*and $P(\delta^2(m)(n-|m|+1)/(n\sigma^2(m)) < \exp(-t))$ is also bounded by the expression on the right-hand side of the preceding display. Clearly, that expression is not larger than $\exp[-((n-|m|)/2)t^2/(t+2)]$.*

PROOF. The case $|m| = 1$ is trivial, as then $\delta^2(m) = \sigma^2(m)$ by Proposition 2.1, and the probabilities of interest reduce to $P(1 > \exp(t))$ and $P(1 < \exp(-t))$, which are both equal to zero. Assume, henceforth, that $|m| > 1$. Let $A$ and $B$ be independent and distributed as $A \sim \chi^2_{|m|-1}$ and $B \sim \chi^2_{n-|m|+1}$. Then, $\delta^2(m)$ is distributed as $\sigma^2(m)(1 + A/B)$ by Proposition 2.1, and

$$\sigma^2(m)(1+A/B)\frac{n-|m|+1}{n\sigma^2(m)} = \frac{|m|-1}{n}\left(\frac{A(n-|m|+1)}{B(|m|-1)} - 1\right) + 1$$

(as is elementary to verify).

First, consider $P(\delta^2(m)(n-|m|+1)/(n\sigma^2(m)) > \exp(t))$. In view of the consideration in the preceding paragraph, this probability equals

$$P\left(\frac{A(n-|m|+1)}{B(|m|-1)} - 1 > \frac{n}{|m|-1}(e^t - 1)\right).$$

Using Lemma A.1 of [20] [with $|m|-1$, $n-|m|+1$ and $(\exp(t)-1)n/(n-|m|+1)$ replacing $a$, $b$ and $\varepsilon$, resp.], the probability in the preceding display is not larger than

$$\exp\left[-\frac{n-|m|+1}{2}\mathcal{K}\left(\frac{|m|-1}{n-|m|+1}, (e^t-1)\frac{n}{n-|m|+1}\right)\right],$$

where the function $\mathcal{K}(r,c)$ is defined for $r > 0$ and $c > -r$ by $\mathcal{K}(r,c) = (1+r)\log((1+r+c)/(1+r)) - r\log((r+c)/r)$. We need to show that the factor involving the $\mathcal{K}$-function in the preceding display satisfies

$$\mathcal{K}\left(\frac{|m|-1}{n-|m|+1}, (e^t-1)\frac{n}{n-|m|+1}\right) \geq \frac{t^2}{t+1+(|m|-1)/n}.$$

To this end, write $s$ as shorthand for $(|m|-1)/n$ and note that we always have $0 < s < 1$. With this, the relation in the preceding display is equivalent to

$$\frac{1}{1-s}\left(t - s\log\frac{e^t+s-1}{s}\right) \geq \frac{t^2}{t+1+s}$$

(expand the formula for the $\mathcal{K}$-function and simplify). It now follows from part (i) of Lemma B.2 that the relation in the preceding display holds.

Now, consider $P(\hat{\delta}^2(m)(n-|m|+1)/(n\sigma^2(m)) < \exp(-t))$, or, equivalently,

(13) $$P\left(\frac{A(n-|m|+1)}{B(|m|-1)} - 1 < \frac{n}{|m|-1}(e^{-t}-1)\right).$$

In case $(e^{-t}-1)n/(|m|-1) \leq -1$, this probability is zero and hence trivially bounded as claimed. In the case where $(e^{-t}-1)n/(|m|-1) > -1$, or, equivalently, $t < -\log(1-s)$, we argue as in the preceding paragraph, mutatis mutandis, to see that (13) is bounded as claimed if

$$\frac{1}{1-s}\left(-t - s\log\frac{e^{-t}+s-1}{s}\right) \geq \frac{t^2}{t+1+s}.$$

But this relation follows by first applying part (ii) and then part (i) of Lemma B.2 as before. □

LEMMA B.4. *For each $t \geq 0$, we have*

$$P\left(\hat{\rho}^2(m)\frac{n-|m|+1}{n\sigma^2(m)} > \exp(t)\right) \leq \exp\left[-\frac{n-|m|}{2}\frac{t^2}{t+2}\right],$$

*and $P(\hat{\rho}^2(m)(n-|m|+1)/(n\sigma^2(m)) < \exp(-t))$ is also bounded by the expression on the right-hand side of the preceding display. The result continues to hold with $\hat{\delta}^2(m)$ replacing $\hat{\rho}^2(m)$.*

PROOF. For $B \sim \chi^2_{n-|m|}$, we have $\hat{\sigma}^2(m) \sim \sigma^2(m)B/(n-|m|)$ (see Proposition 2.1). Hence, $\hat{\rho}^2(m) = n\hat{\sigma}^2(m)/(n-|m|+1)$ is distributed as $(\sigma^2(m)B/(n-|m|))n/(n-|m|+1)$, and

$$\frac{\sigma^2(m)B}{n-|m|}\frac{n}{n-|m|+1}\frac{n-|m|+1}{n\sigma^2(m)} = \frac{B}{n-|m|}.$$

First consider $P(\hat{\rho}^2(m)(n-|m|+1)/(n\sigma^2(m)) > \exp(t))$. By the preceding consideration, this probability equals

$$P\left(\frac{B}{n-|m|} - 1 > \exp(t) - 1\right).$$

Using Lemma A.2 of [20] [with $n - |m|$ and $\exp(t) - 1$ replacing $b$ and $\varepsilon$, resp.], we see that the probability in the preceding display is not larger than

$$\exp\left[-\frac{n-|m|}{2}(e^t - 1 - t)\right].$$

Now, Lemma B.2(iii) entails that the expression in the preceding display is bounded by $\exp[-(n-|m|)t^2/(2(t+2))]$ as required. The derivation of the upper bound for $P(\hat{\rho}^2(m)(n-|m|+1)/(n\sigma^2(m)) < \exp(-t))$ is completely analogous. Finally, the statement in parentheses follows, because $\hat{\rho}^2(m)$ and $\hat{\delta}^2(m)$ are given by the same formula. $\square$

LEMMA B.5. *For each $t \geq 0$, we have*

$$P\left(\rho(m)^2 \frac{n-|m|+1}{n\sigma^2(m)} < \exp(-t)\right) \leq \exp\left[-\frac{n-|m|}{2}\frac{t^2}{t+2}\right]$$

*and*

$$P\left(\rho(m)^2 \frac{n-|m|+1}{n\sigma^2(m)} > \exp(t)\right) \leq 3\exp\left[-\frac{n-|m|}{4}\frac{t^2}{t+4}\right].$$

PROOF. The first inequality follows immediately from Lemma B.3 upon noting that $\delta^2(m) \leq \delta^2(m) + \nu^2(m) = \rho(m)^2$.

The second inequality holds trivially in case the upper bound is larger than one. We exclude the trivial case and hence assume that $\log(3) < ((n-|m|)/4)t^2/(t+4)$. For later use, we note that this entails that $1 < nt/2$ [because $\log(3) < (n/4)t^2/(t+4) < (n/4)t$, so that $2\log(3) < nt/2$, where the lower bound is larger than one]. In the second inequality of the lemma, the expression on the left-hand side is bounded by

(14)
$$P\left(\delta^2(m)\frac{n-|m|+1}{n\sigma^2(m)} > \exp(t/2)\right)$$
$$+ P\left(\nu^2(m)\frac{n-|m|+1}{n\sigma^2(m)} > \frac{t}{2}\exp(t/2)\right),$$

because $\rho^2(m) = \nu^2(m) + \delta^2(m)$ and $e^t = e^{t/2}e^{t/2} \geq e^{t/2}(1 + t/2)$. The first term in (14) is bounded by $\exp[-((n-|m|)/4)t^2/(t+4)]$ (use Lemma B.3 with $t/2$ replacing $t$ and simplify).

To complete the proof, we need to show that the second term in (14) is bounded by $2\exp[-((n-|m|)/4)t^2/(t+4)]$. To this end, recall from Proposition 2.1 that $\nu^2(m)$ is distributed as $(A/n)\delta^2(m)$, where $A \sim \chi_1^2$ independent of $\delta^2(m)$. Hence, the second term in (14) is bounded by

$$P\left(A > n\frac{t}{2}\right) + P\left(\delta^2(m)\frac{n-|m|+1}{n\sigma^2(m)} > \exp(t/2)\right).$$

The second term in the preceding display coincides with the first term in (14) and is bounded by $\exp[-((n-|m|)/4)t^2/(t+4)]$ as shown before. To complete the proof, we need to show that the first term is also bounded by that quantity. By the second inequality of Lemma B.1, the term in question is bounded by $\sqrt{2/\pi}\exp[-nt/4 - \log(nt/2)/2]$. Now recall that we have $nt/2 > 1$ and note that $\sqrt{2/\pi} < 1$. Hence, the first term in the preceding display is bounded by $\exp[-nt/4] \leq \exp[-((n-|m|)/4)t^2/(t+4)]$. □

## APPENDIX C: PROOFS FOR SECTION 3

PROOF OF THEOREM 3.1. We first derive separate upper bounds for $P(\rho^2(m)/\hat{\rho}^2(m) < e^{-\varepsilon})$ and for $P(\rho^2(m)/\hat{\rho}^2(m) > e^\varepsilon)$.

If $\rho^2(m)/\hat{\rho}^2(m) < e^{-\varepsilon}$, then either

$$\rho^2(m)\frac{n-|m|+1}{n\sigma^2(m)} < \exp(-\varepsilon/2) \quad \text{or} \quad \hat{\rho}^2(m)\frac{n-|m|+1}{n\sigma^2(m)} > \exp(\varepsilon/2).$$

Using Lemma B.5 to bound the probability of the first event in the preceding display and using Lemma B.4 to bound the probability of the second one, we see that

$$P\left(\frac{\rho^2(m)}{\hat{\rho}^2(m)} < e^{-\varepsilon}\right) < 2\exp\left[-\frac{n-|m|}{4}\frac{\varepsilon^2}{\varepsilon+4}\right].$$

Clearly, the upper bound in the preceding display is not larger than $2\exp[-((n-|m|)/8)\varepsilon^2/(\varepsilon+8)]$.

For $P(\rho^2(m)/\hat{\rho}^2(m) > e^\varepsilon)$, we argue similarly as in the preceding paragraph to obtain

$$P\left(\frac{\rho^2(m)}{\hat{\rho}^2(m)} > e^\varepsilon\right) < 4\exp\left[-\frac{n-|m|}{8}\frac{\varepsilon^2}{\varepsilon+8}\right].$$

Relation (4) follows from this. □

COROLLARY C.1. *In the setting of Theorem 3.1, relation (4) continues to hold with $\delta^2(m)$ replacing $\rho^2(m)$; in that case, the constants 6 and 8 in (4) can both be replaced by 4.*

PROOF. The result follows by arguing as in the proof of Theorem 3.1, mutatis mutandis, now using Lemma B.3 instead of Lemma B.5. □

PROOF OF COROLLARY 3.2. Let $E$ denote the event where

$$\max_{m \in \mathcal{M}} \left| \log \frac{\hat{\rho}^2(m)}{\rho^2(m)} \right| \leq \varepsilon/2,$$

and note that the complement of $E$, that is, $E^c$, is such that

$$P(E^c) \leq \sum_{m \in \mathcal{M}} P\left( \left| \log \frac{\hat{\rho}^2(m)}{\rho^2(m)} \right| > \varepsilon/2 \right) \leq \sum_{m \in \mathcal{M}} 6 \exp\left[ -\frac{n-|m|}{8} \frac{(\varepsilon/2)^2}{(\varepsilon/2)+8} \right]$$

$$\leq 6\#\mathcal{M} \exp\left[ -\frac{n-|\mathcal{M}|}{16} \frac{\varepsilon^2}{\varepsilon+16} \right].$$

(In the preceding chain of inequalities, the first one is derived from Bonferroni's inequality, the second one follows from Theorem 3.1, and the last one is obvious in view of the definitions of $\#\mathcal{M}$ and $|\mathcal{M}|$.)

To derive the first statement of the corollary, first note that the relation $0 \leq \log(\rho^2(\hat{m})/\rho^2(m_\rho))$ is always satisfied. Moreover, observe that

$$\log \frac{\rho^2(\hat{m})}{\rho^2(m_\rho)} = \log \frac{\rho^2(\hat{m})}{\hat{\rho}^2(\hat{m})} + \log \frac{\hat{\rho}^2(\hat{m})}{\hat{\rho}^2(m_\rho)} + \log \frac{\hat{\rho}^2(m_\rho)}{\rho^2(m_\rho)}$$

almost surely, because the event where $\hat{\rho}^2(m) > 0$ for each $m \in \mathcal{M}$ has probability one. On the right-hand side of the preceding equality, the second term is nonpositive, and hence

$$\log \frac{\rho^2(\hat{m})}{\rho^2(m_\rho)} \leq \left| \log \frac{\rho^2(\hat{m})}{\hat{\rho}^2(\hat{m})} \right| + \left| \log \frac{\hat{\rho}^2(m_\rho)}{\rho^2(m_\rho)} \right|$$

almost surely. Hence, on the event $E$, we see that $\log(\rho^2(\hat{m})/\rho^2(m_\rho))$ is between zero and $\varepsilon$, and $P(E^c)$ is bounded from above as required.

For the second statement of the corollary, define $E$ as before but now with $\varepsilon$ replacing $\varepsilon/2$. It is easy to see that now $P(E^c)$ is not larger than $6\#\mathcal{M}\exp[-((n-|\mathcal{M}|)/8)\varepsilon^2/(\varepsilon+8)]$. On the event $E$, we clearly have $|\log \hat{\rho}^2(\hat{m})/\rho^2(\hat{m})| \leq \varepsilon$. □

## APPENDIX D: PROOFS FOR SECTION 4

The following lemma provides an upper bound for the total variation distance of two normal distributions in terms of their parameters and will be instrumental in the proof of Theorem 4.1. We believe that the lemma is well known, in some form or another, but we could not find an appropriate reference.

LEMMA D.1. *Write $N(a, s^2)$ and $N(0,1)$ for the Gaussian measures with the indicated parameters (where $a \in \mathbb{R}$, $s > 0$). Then the total variation distance of these two measures is bounded as*

$$\|N(a, s^2) - N(0, 1)\|_{\mathrm{TV}} \leq \frac{|a|}{\sqrt{2\pi}} + \frac{|\log(s^2)|}{\sqrt{2\pi e}}.$$

REMARK D.1. Of course $\|N(a, s^2) - N(0, 1)\|_{\mathrm{TV}}$ is trivially bounded by one. Moreover, the lemma also entails that that total variation distance is also bounded by $|a/s|/\sqrt{2\pi} + |\log(s^2)|/\sqrt{2\pi e}$, because $\|N(a, s^2) - N(0, 1)\|_{\mathrm{TV}} = \|N(0, 1) - N(-a/s, 1/s^2)\|_{\mathrm{TV}}$.

PROOF OF LEMMA D.1. Recall that the total variation distance of two mutually absolutely continuous probability measures $P$ and $Q$ is given by

$$(15) \qquad \|P - Q\|_{\mathrm{TV}} = P(\log(p/q) > 0) - Q(\log(p/q) > 0),$$

where $p$ and $q$ are the densities of $P$ and $Q$, respectively, with respect to a common dominating sigma-finite measure. Write $\phi(t)$ for the Lebesgue density of $N(0, 1)$, and note that the Lebesgue density of $N(a, s^2)$ is then given by $\phi((t-a)/s)/s$. The log-likelihood ratio of $N(a, s^2)$ and $N(0, 1)$ in hence given by

$$(16) \qquad \log\left(\frac{\phi((t-a)/s)/s}{\phi(t)}\right) = \frac{1}{2}\log(1/s^2) - \frac{1}{2}\left(\frac{(t-a)^2}{s^2} - t^2\right).$$

The total variation distance of $N(a, s^2)$ and $N(0, 1)$ is bounded from above by $\|N(a, s^2) - N(a, 1)\|_{\mathrm{TV}} + \|N(a, 1) - N(0, 1)\|_{\mathrm{TV}}$ or, equivalently, by

$$(17) \qquad \|N(a, 1) - N(0, 1)\|_{\mathrm{TV}} + \|N(0, s^2) - N(0, 1)\|_{\mathrm{TV}}.$$

The proof will be complete if we can show that the first term in (17) is bounded by $|a|/\sqrt{2\pi}$ and that the second term in (17) is bounded by $|\log(s^2)|/\sqrt{2\pi e}$.

To bound the first term in (17), we first use (16) with $s^2$ replaced by 1 to see that the log-likelihood ratio of $N(a, 1)$ and $N(0, 1)$ is given by $(-1/2)(-2ta + a^2)$, which is positive if and only if $ta > a^2/2$. Using (15) with $N(a, 1)$ and $N(0, 1)$ replacing $P$ and $Q$, respectively, it is elementary to verify that

$$\|N(a, 1) - N(0, 1)\|_{\mathrm{TV}} = 2\Phi(|a|/2) - 1,$$

where $\Phi(\cdot)$ denotes the standard Gaussian c.d.f. For $x \geq 0$, set $f(x) = 2\Phi(x) - 1$ and $g(x) = x\sqrt{2/\pi}$. If we can show that $f(x) \leq g(x)$, it will follow that the expression in the preceding display is bounded from above by $g(|a|/2) = |a|/\sqrt{2\pi}$. To show that $f(x) \le g(x)$ for $x \ge 0$, note that $f(0) = g(0) = 0$, and that

$$f'(x) = 2\frac{1}{\sqrt{2\pi}} e^{-x^2/2} \le \sqrt{2/\pi} = g'(x).$$

The claim now follows because $f(x) = f(0) + \int_0^x f(t)\,dt \le g(0) + \int_0^x g(t)\,dt = g(x)$.

For the second term in (17), note that it suffices to consider the case where $s^2 > 1$ [because that term is trivially bounded from above by $|\log(s^2)|/\sqrt{2\pi e}$ in case $s^2 = 1$; because $\|N(0,s^2) - N(0,1)\|_{TV} = \|N(0,1) - N(0,1/s^2)\|_{TV}$; and because $|\log(s^2)| = |\log(1/s^2)|$]. Use (16) with $a$ replaced by 0 to see that the log-likelihood ratio of $N(0,s^2)$ and $N(0,1)$ is given by $-\log(s^2)/2 - t^2(1/s^2 - 1)/2$. This log-likelihood ratio is positive at $t$ if and only if $t^2 > s^2 \log(s^2)/(s^2 - 1)$, because $s^2 > 1$. Using (15) with $N(0,s^2)$ and $N(0,1)$ replacing $P$ and $Q$, respectively, it is now easy to see that

$$\|N(0,s^2) - N(0,1)\|_{TV} = F\left(s^2 \frac{\log(s^2)}{s^2 - 1}\right) - F\left(\frac{\log(s^2)}{s^2 - 1}\right),$$

where $F(\cdot)$ denotes the c.d.f. of a chi-square distributed random variable with one degree of freedom. Using the first inequality of Lemma B.1 with $s^2$ replacing $t$, we see that the expression in the preceding display is bounded by $\log(s^2)/\sqrt{2\pi e}$ as required. $\square$

PROOF OF THEOREM 4.1. Because $\|N(0,\hat{\delta}^2(m)) - N(\nu(m),\delta^2(m))\|_{TV} = \|N(-\nu(m)/\delta(m), \hat{\delta}^2(m)/\delta^2(m)) - N(0,1)\|_{TV}$, Lemma D.1 entails that

$$\|N(0,\hat{\delta}^2(m)) - N(\nu(m),\delta^2(m))\|_{TV}$$
$$\le \frac{|\nu(m)/\delta(m)|}{\sqrt{2\pi}} + \frac{|\log(\hat{\delta}^2(m)/\delta^2(m))|}{\sqrt{2\pi e}}.$$

In view of this, and because $1/\sqrt{2\pi e} < 1/4$, we get

(18)
$$P\left(\|N(0,\hat{\delta}^2(m)) - N(\nu(m),\delta^2(m))\|_{TV} > \frac{1}{\sqrt{n}} + \varepsilon\right)$$
$$\le P\left(\frac{|\nu(m)/\delta(m)|}{\sqrt{2\pi}} + \frac{|\log(\hat{\delta}^2(m)/\delta^2(m))|}{\sqrt{2\pi e}} > \frac{1}{\sqrt{n}} + \varepsilon\right)$$
$$\le P\left(\frac{|\nu(m)/\delta(m)|}{\sqrt{2\pi}} > \frac{1}{\sqrt{n}} + \frac{\varepsilon}{2}\right) + P\left(\frac{|\log(\hat{\delta}^2(m)/\delta^2(m))|}{\sqrt{2\pi e}} > \frac{\varepsilon}{2}\right)$$
$$\le P\left(\left|\frac{\nu(m)}{\delta(m)}\right| > \sqrt{2\pi}\left(\frac{1}{\sqrt{n}} + \frac{\varepsilon}{2}\right)\right) + P\left(\left|\log\left(\frac{\hat{\delta}^2(m)}{\delta^2(m)}\right)\right| > 2\varepsilon\right).$$

For the second term in (18), recall that $\hat{\hat{\delta}}^2(m) = \hat{\rho}^2(m)$ and use Corollary C.1 to obtain

$$P\left(\left|\log\left(\frac{\hat{\hat{\delta}}^2(m)}{\delta^2(m)}\right)\right| > 2\varepsilon\right) \leq 4\exp\left[-\frac{n-|m|}{2}\frac{\varepsilon^2}{\varepsilon+2}\right].$$

To complete the proof, we need to show that the first term in (18) is bounded by $3\exp[-((n-|m|)/2)\varepsilon^2/(\varepsilon+2)]$. For the first term in (18), observe that

$$P\left(\frac{\nu^2(m)}{\delta^2(m)} > 2\pi\left(\frac{1}{\sqrt{n}} + \frac{\varepsilon}{2}\right)^2\right)$$

$$\leq P\left(\nu^2(m)\frac{n-|m|+1}{n\sigma^2(m)} > 2\pi\left(\frac{1}{\sqrt{n}} + \frac{\varepsilon}{2}\right)^2 e^{-\varepsilon}\right)$$

$$+ P\left(\frac{1}{\delta^2(m)}\frac{n\sigma^2(m)}{n-|m|+1} > e^{\varepsilon}\right).$$

In the preceding display, the second term on the right-hand side equals $P(\delta^2(m) \times (n-|m|+1)/(n\sigma^2(m)) < e^{-\varepsilon}) \leq \exp[-\frac{n-|m|}{2}\frac{\varepsilon^2}{\varepsilon+2}]$, where inequality follows from Lemma B.3. It remains to show that, in the preceding display, the first term on the right is not larger than $2\exp[-((n-|m|)/2)\varepsilon^2/(\varepsilon+2)]$. To this end, let $A \sim \chi_1^2$ independent of $\delta^2(m)$. In view of Proposition 2.1, $\nu^2(m)$ has the same distribution as $(A/n)\delta^2(m)$. Therefore,

$$P\left(\nu^2(m)\frac{n-|m|+1}{n\sigma^2(m)} > 2\pi\left(\frac{1}{\sqrt{n}} + \frac{\varepsilon}{2}\right)^2 e^{-\varepsilon}\right)$$

(19)
$$= P\left(\frac{A}{n}\delta^2(m)\frac{n-|m|+1}{n\sigma^2(m)} > 2\pi\left(\frac{1}{\sqrt{n}} + \frac{\varepsilon}{2}\right)^2 e^{-\varepsilon}\right)$$

$$\leq P\left(\frac{A}{n} > 2\pi\left(\frac{1}{\sqrt{n}} + \frac{\varepsilon}{2}\right)^2 e^{-2\varepsilon}\right) + P\left(\delta^2(m)\frac{n-|m|+1}{n\sigma^2(m)} > e^{\varepsilon}\right).$$

For the second term on the far right-hand side of (19), we again use Lemma B.3 to get

$$P\left(\delta^2(m)\frac{n-|m|+1}{n\sigma^2(m)} > e^{\varepsilon}\right) \leq \exp\left[-\frac{n-|m|}{2}\frac{\varepsilon^2}{\varepsilon+2}\right].$$

In view of this, the proof will be complete if the first term on the far right of (19) is bounded by $\exp[-((n-|m|)/2)\varepsilon^2/(\varepsilon+2)]$.

For the first term on the far right of (19), we have

(20) $$P\left(\frac{A}{n} > 2\pi\left(\frac{1}{\sqrt{n}} + \frac{\varepsilon}{2}\right)^2 e^{-2\varepsilon}\right) \leq \exp\left[-\pi n\left(\frac{1}{\sqrt{n}} + \frac{\varepsilon}{2}\right)^2 e^{-2\varepsilon}\right],$$

in view of the second inequality of Lemma B.1 and because $2\pi n(1/\sqrt{n} + \varepsilon/2)^2 e^{-2\varepsilon} = 2\pi(1 + \sqrt{n}\varepsilon/2)^2 e^{-2\varepsilon} \geq 2\pi e^{-2\varepsilon} \geq 2\pi e^{-2\log 2} = 2\pi/4 > 1$. We now show that the right-hand side of (20) is not larger than $\exp[-(n/2)\varepsilon^2/(\varepsilon+2)]$ or, equivalently, that

$$(21) \qquad 2\pi\left(\frac{1}{\sqrt{n}} + \frac{\varepsilon}{2}\right)^2 e^{-2\varepsilon} \geq \frac{\varepsilon^2}{\varepsilon+2}.$$

In this inequality, the left-hand side satisfies

$$2\pi(1/\sqrt{n} + \varepsilon/2)^2 e^{-2\varepsilon} > (\pi/2)\varepsilon^2 e^{-2\varepsilon}.$$

Thus, (21) will follow if $\pi(\varepsilon+2) \geq 2e^{2\varepsilon}$, or if $f(\varepsilon) = \pi(\varepsilon+2) - 2e^{2\varepsilon} \geq 0$. It is now easy to verify that $f(\varepsilon)$ is strictly decreasing and that $f(\log(2)) > 0$. Hence, the expression on the left of (20) or, equivalently, the first term on the far right of (19) is bounded by $\exp[-(n/2)\varepsilon^2/(\varepsilon+2)] < \exp[-((n-|m|)/2)\varepsilon^2/(\varepsilon+2)]$. $\square$

PROOF OF COROLLARY 4.2. Using Bonferroni's inequality and Theorem 4.1, this result follows immediately by arguing as in the first paragraph of the proof of Corollary 3.2. $\square$

PROOF OF PROPOSITION 4.3. For measurable $A \subseteq \mathbb{R}$, write $\mathbb{L}(\hat{m}; A)$ and $\hat{\mathbb{L}}(\hat{m}; A)$ for the probability of $A$ under $\mathbb{L}(\hat{m})$ and under $\hat{\mathbb{L}}(\hat{m})$, respectively. We have $y^{(f)} \in \mathcal{I}(\hat{m})$ if and only if $\hat{y}^{(f)}(\hat{m}) - y^{(f)}$ lies in the interval $[-q_\alpha \hat{\delta}(\hat{m}), q_\alpha \hat{\delta}(\hat{m})]$. Writing $A$ as shorthand for that interval, the conditional coverage probability of $\mathcal{I}(\hat{m})$ equals $\mathbb{L}(\hat{m}; A)$. Because the interval $\mathcal{I}(\hat{m})$ is constructed with nominal coverage probability $1 - \alpha$ assuming that $\hat{y}^{(f)}(\hat{m}) - y^{(f)}$ is distributed as $\hat{\mathbb{L}}(\hat{m})$, we have $\hat{\mathbb{L}}(\hat{m}; A) = 1 - \alpha$. The result now follows immediately from Corollary 4.2. $\square$

PROOF OF PROPOSITION 4.4. To bound $P(\log \hat{\delta}(\hat{m})/\delta(m_\delta) > \varepsilon)$, we first note that

$$\log \frac{\hat{\delta}^2(\hat{m})}{\delta^2(m_\delta)} = \log \frac{\hat{\delta}^2(\hat{m})}{\hat{\delta}^2(m_\delta)} + \log \frac{\hat{\delta}^2(m_\delta)}{\delta^2(m_\delta)} \leq \log \frac{\hat{\delta}^2(m_\delta)}{\delta^2(m_\delta)}$$

almost surely, because $\hat{m}$ is a minimizer of $\hat{\delta}^2(\cdot) = \hat{\rho}^2(\cdot)$, and because the event where $\hat{\delta}^2(m) > 0$ for each $m \in \mathcal{M}$ has probability one. Hence, $P(\log \hat{\delta}(\hat{m})/\delta(m_\delta) > \varepsilon)$ or, equivalently, $P(\log \hat{\delta}^2(\hat{m})/\delta^2(m_\delta) > 2\varepsilon)$, is bounded by

$$\sum_{m \in \mathcal{M}} P\left(\log \frac{\hat{\delta}^2(m)}{\delta^2(m)} > 2\varepsilon\right)$$

in view of Bonferroni's inequality.

Similarly, to bound $P(\log \hat{\delta}(\hat{m})/\delta(m_\delta) < -\varepsilon)$, we observe that

$$\log \frac{\hat{\delta}^2(\hat{m})}{\delta^2(m_\delta)} = \log \frac{\hat{\delta}^2(\hat{m})}{\delta^2(\hat{m})} + \log \frac{\delta^2(\hat{m})}{\delta^2(m_\delta)} \geq \log \frac{\hat{\delta}^2(\hat{m})}{\delta^2(\hat{m})},$$

because $m_\delta$ is a minimizer of $\delta^2(\cdot)$. Arguing similarly as in the preceding paragraph, we thus see that $P(\log \hat{\delta}^2(\hat{m})/\delta^2(m_\delta) < -2\varepsilon)$ is bounded from above by

$$\sum_{m \in \mathcal{M}} P\bigg(\log \frac{\hat{\delta}^2(m)}{\delta^2(m)} < -2\varepsilon\bigg).$$

Adding the bounds for $P(\log \hat{\delta}(\hat{m})/\delta(m_\delta) > \varepsilon)$ and for $P(\log \hat{\delta}(\hat{m})/\delta(m_\delta) < -\varepsilon)$ obtained so far, we see that $P(|\log \hat{\delta}(\hat{m})/\delta(m_\delta)| > \varepsilon)$ is bounded by

$$\sum_{m \in \mathcal{M}} P\bigg(\bigg|\log \frac{\hat{\delta}^2(m)}{\delta^2(m)}\bigg| > 2\varepsilon\bigg).$$

Recalling that $\hat{\delta}^2(\cdot) = \hat{\rho}^2(\cdot)$, the result now follows from Corollary C.1. $\square$

**Acknowledgments.** I would like to thank John Hartigan, Richard Nickl, David Pollard and Benedikt Pötscher for inspiring critique and helpful feedback.

## REFERENCES

[1] ADAM, B.-L., QU, Y., DAVIS, J. W., WARD, M. D., CLEMENTS, M. A., CAZARES, L. H., SEMMES, O. J., SCHELLMANNER, P. F., YASUI, Y., FENG, Z. and WRIGHT, G. L. J. (2002). Serum protein fingerprinting coupled with a pattern-matching algorithm distinguishes prostate cancer from benign prostate hyperplasia and healthy men. *Cancer Research* **62** 3609–3614.
[2] BARAUD, Y. (2004). Confidence balls in Gaussian regression. *Ann. Statist.* **32** 528–551. MR2060168
[3] BARNDORFF-NIELSEN, O. E. and COX, D. R. (1996). Prediction and asymptotics. *Bernoulli* **2** 319–340. MR1440272
[4] BERAN, R. and DÜMBGEN, L. (1998). Modulation of estimators and confidence sets. *Ann. Statist.* **26** 1826–1856. MR1673280
[5] BREIMAN, L. and FREEDMAN, D. (1983). How many variables should be entered in a regression equation? *J. Amer. Statist. Assoc.* **78** 131–136. MR0696857
[6] CAI, T. T. and LOW, M. G. (2004). An adaptation theory for nonparametric confidence intervals. *Ann. Statist.* **32** 1805–1840. MR2102494
[7] CAI, T. T. and LOW, M. G. (2006). Adaptive confidence balls. *Ann. Statist.* **34** 202–228. MR2275240
[8] DING, A. A. and HWANG, J. T. G. (1999). Prediction intervals, factor analysis models, and high-dimensional empirical linear prediction. *J. Amer. Statist. Assoc.* **94** 446–455. MR1702316
[9] GEISSER, S. (1993). *Predictive Inference: An Introduction. Monographs on Statistics and Applied Probability* **55**. Chapman & Hall, New York. MR1252174
[10] GENOVESE, C. R. and WASSERMAN, L. (2005). Confidence sets for nonparametric wavelet regression. *Ann. Statist.* **33** 698–729. MR2163157
[11] GENOVESE, C. R. and WASSERMAN, L. (2008). Adaptive confidence bands. *Ann. Statist.* **36** 875–905. MR2396818

[12] Golub, T. R., Slonim, D. K., Tamayo, P., Huard, C., Gaasenbeek, M., Mesirov, J. P., Coller, H., Loh, M. L., Downing, J. R., Caligiuri, M. A., Bloomfield, D. C. and Lander, E. S. (1999). Molecular classification of cancer: Class discovery and class prediction by gene expression monitoring. *Science* **286** 531–537.

[13] Hocking, R. R. (1976). The analysis and selection of variables in linear regression. *Biometrics* **32** 1–49. MR0398008

[14] Hoffmann, M. and Lepski, O. (2002). Random rates in anisotropic regression. *Ann. Statist.* **30** 325–396. MR1902892

[15] Joshi, V. M. (1969). Admissibility of the usual confidence sets for the mean of a univariate or bivariate normal population. *Ann. Math. Statist.* **40** 1042–1067. MR0264811

[16] Juditsky, A. and Lambert-Lacroix, S. (2003). Nonparametric confidence set estimation. *Math. Methods Statist.* **12** 410–428. MR2054156

[17] Kabaila, P. and Leeb, H. (2006). On the large-sample minimal coverage probability of confidence intervals after model selection. *J. Amer. Statist. Assoc.* **101** 619–629. MR2256178

[18] Leeb, H. (2005). The distribution of a linear predictor after model selection: Conditional finite-sample distributions and asymptotic approximations. *J. Statist. Plann. Inference* **134** 64–89. MR2146086

[19] Leeb, H. (2006). The distribution of a linear predictor after model selection: Unconditional finite-sample distributions and asymptotic approximations. *IMS Lecture Notes—Monograph Series* **49** 291–311. MR2338549

[20] Leeb, H. (2008). Evaluation and selection of models for out-of-sample prediction when the sample size is small relative to the complexity of the data-generating process. *Bernoulli* **14** 661–690.

[21] Leeb, H. and Pötscher, B. M. (2003). The finite-sample distribution of post-model-selection estimators, and uniform versus non-uniform approximations. *Econometric Theory* **19** 100–142. MR1965844

[22] Leeb, H. and Pötscher, B. M. (2005). Model selection and inference: Facts and fiction. *Econometric Theory* **21** 21–59. MR2153856

[23] Leeb, H. and Pötscher, B. M. (2006). Can one estimate the conditional distribution of post-model-selection estimators? *Ann. Statist.* **34** 2554–2591. MR2291510

[24] Leeb, H. and Pötscher, B. M. (2008). Can one estimate the unconditional distribution of post-model-selection estimators? *Econometric Theory* **24** 338–376.

[25] Li, K.-C. (1989). Honest confidence regions for nonparametric regression. *Ann. Statist.* **17** 1001–1008. MR1015135

[26] Nychka, D. (1988). Bayesian confidence intervals for smoothing splines. *J. Amer. Statist. Assoc.* **83** 1134–1143. MR0997592

[27] Pötscher, B. M. (1991). Effects of model selection on inference. *Econometric Theory* **7** 163–185. MR1128410

[28] Robins, J. and van der Vaart, A. (2006). Adaptive nonparametric confidence sets. *Ann. Statist.* **34** 229–253. MR2275241

[29] Shen, X., Huang, H.-C. and Ye, J. (2004). Inference after model selection. *J. Amer. Statist. Assoc.* **99** 751–761. MR2090908

[30] Souders, T. M. and Stenbakken, G. N. (1991). Cutting the high cost of testing. *IEEE Spectrum* **28** 48–51.

[31] Stenbakken, G. N. and Souders, T. M. (1987). Test point selection and testability measures via QR factorization of linear models. *IEEE Trans. Instrum. Meas.* **36** 406–410.

[32] THOMPSON, M. L. (1978). Selection of variables in multiple regression: Part II. Chosen procedures, computations and examples. *Int. Statist. Rev.* **46** 129–146. MR0514059
[33] TIBSHIRANI, R., SAUNDERS, M., ROSSET, S., ZHU, J. and KNIGHT, K. (2005). Sparsity and smoothness via the fused lasso. *J. Roy. Statist. Soc. Ser. B* **67** 91–108. MR2136641
[34] VAN DE VIJVER, M. J., HE, Y. D., VAN'T VEER, L. J., DAI, H., HART, A. A. M., VOSKUIL, D. W., SCHREIBER, G. J., PETERSE, J. L., ROBERTS, C., MARTON, M. J., PARRISH, M., ATSMA, D., WITTEVEEN, A., GLAS, A., DELAHAYE, L., VAN DER VELDE, T., BARTELINK, H., RODENHUIS, S., RUTGERS, E. T., FRIEND, S. H. and BERNARDS, R. (2002). A gene-expression signature as a predictor of survival in breast cancer. *The New England Journal of Medicine* **347** 1999–2009.
[35] VAN'T VEER, L. J., DAI, H., VAN DE VIJVER, M. J., HE, Y. D., HART, A. A. M., MAO, M., PETERSE, H. L., VAN DER KOOY, K., MARTON, M. J., WITTEVEEN, A. T., SCHREIBER, G. J., KERKHOVEN, R. M., ROBERTS, C., LINSLEY, P. S., BERNARDS, R. and FRIEND, S. H. (2002). Gene expression profiling predicts clinical outcome of breast cancer. *Nature* **415** 530–536.
[36] WAHBA, G. (1983). Bayesian "confidence intervals" for the cross-validated smoothing spline. *J. Amer. Statist. Assoc.* **45** 133–150. MR0701084
[37] WEST, M., BLANCHETTE, C., DRESSMAN, H., HUANG, E., ISHIDA, S., SPANG, R., ZUZAN, H., OLSON, J. A. J., MARKS, J. R. and NEVINS, J. R. (2001). Predicting the clinical status of human breast cancer by using gene expression profiles. *Proc. Natl. Acad. Sci. U.S.A.* **98** 11462–11467.

DEPARTMENT OF STATISTICS
YALE UNIVERSITY
24 HILLHOUSE AVENUE
NEW HAVEN, CONNECTICUT 06511
USA
E-MAIL: hannes.leeb@yale.edu